\documentclass[AMA,STIX1COL]{WileyNJD-v2}
\usepackage{moreverb}

\articletype{Article Type}%

\received{<day> <Month>, <year>}
\revised{<day> <Month>, <year>}
\accepted{<day> <Month>, <year>}

\usepackage{amsmath,amssymb,amsfonts}
\usepackage{graphicx}
\usepackage{textcomp}

\usepackage{graphicx}      

\usepackage{subcaption}
\usepackage{mwe}

\usepackage{booktabs}
\usepackage{moreverb}

\usepackage{amsmath} 
\usepackage{amssymb}  

\usepackage{mathtools, cuted}

\newcommand{\dd}{{\mathrm{d}}}
\newcommand{\barr}{\begin{array}}
	\newcommand{\earr}{\end{array}}
\newcommand{\bvec}{ \left[ \!\! \barr{cccccccccccc} }
\newcommand{\evec}{ \earr \!\! \right] }

\newcommand{\R}{{\mathbb{R}}}
\newcommand{\LL}{{\mathcal L}}

\newcommand{\eye}{\mathbb{1}}
\newcommand{\zero}{\mathbb{0}}

\newcommand{\Ix}{I_\x}

\newcommand{\nm}{n_{\mathrm{m}}}

\newcommand{\software}[1]{\texttt{#1}}
\renewcommand{\min}{\text{min}}
\newcommand{\argmin}{\text{arg min}}
\renewcommand{\st}{\text{s.t.}}
\newcommand{\bx}{\hat{x}_0}

\newcommand{\x}{\mathrm{x}}
\renewcommand{\u}{\mathrm{u}}
\newcommand{\nx}{n_{\mathrm{x}}}
\newcommand{\nck}{n_{\mathrm{c},k}}
\newcommand{\nU}{n_{\mathrm{u}}}

\newcommand{\A}{\mathcal{A}}
\newcommand{\Z}{\mathcal{Z}}
\renewcommand{\H}{\mathcal{H}}
\newcommand{\Hp}{\tilde{\mathcal{H}}}
\newcommand{\V}{\Gamma}
\newcommand{\T}{\mathcal{T}}
\newcommand{\Ta}{\mathcal{T}_{\mathrm{a}}}
\newcommand{\Ts}{\mathcal{T}_{\mathrm{s}}}
\newcommand{\Tc}{\mathcal{T}_{\mathrm{c}}}
\newcommand{\K}{\mathcal{K}}

\newcommand{\La}{L_{\mathrm{a}}}
\newcommand{\Ls}{L_{\mathrm{s}}}

\renewcommand{\P}{\mathcal{P}}

\newcommand{\W}{\mathcal{W}}

\newcommand{\Rho}{P}

\newcommand{\polpar}{d}

\newcommand{\s}{s}
\renewcommand{\S}{S}

\newcommand{\D}{D}

\newcommand{\Tsamp}{T_{\mathrm{s}}}

\renewcommand{\O}{\mathcal{O}}

\usepackage{algorithm}
\usepackage{algpseudocode}
\algtext*{EndFor}
\algtext*{EndIf}
\algrenewcommand\algorithmicrequire{\textbf{Input:}}
\algrenewcommand\algorithmicensure{\textbf{Output:}}

\usepackage{graphicx}

\usepackage{bbold}

\usepackage{color}

\definecolor{wheat}{rgb}{0.96,0.87,0.70}

\newcommand{\change}[1]{#1}

\begin{document}
	
	\title{\software{PRESAS}: Block-Structured Preconditioning of Iterative Solvers within a Primal Active-Set Method for fast MPC\protect\thanks{PRESAS: Preconditioned RESidual method within Active-Set strategy for solving constrained linear-quadratic optimal control problems.}}
	
	\author[1]{Rien Quirynen*}
	
	\author[1]{Stefano Di Cairano}
	
	\authormark{R. QUIRYNEN \textsc{et al}}

	\address[1]{\orgdiv{Control and Dynamical Systems}, \orgname{Mitsubishi Electric Research Laboratories}, \orgaddress{\state{Massachusetts}, \country{USA}}}
	
	
	\corres{*\email{quirynen@merl.com}}
	
	\presentaddress{201 Broadway, 8th Floor, Cambridge, MA 02139-1955}
	
	\abstract[Abstract]{Model predictive control~(MPC) for linear dynamical systems requires solving an optimal control structured quadratic program~(QP) at each sampling instant. This paper proposes a primal active-set strategy\change{, called~\software{PRESAS},} for the efficient solution of such block-sparse QPs, based on a preconditioned iterative solver to compute the search direction in each iteration. Rank-one factorization updates of the preconditioner result in a per-iteration computational complexity of $\O(N m^2)$, where $m$ denotes the number of state and control variables and $N$ the number of control intervals. Three different block-structured preconditioning techniques are presented and their numerical properties are studied further. In addition, an augmented Lagrangian based implementation is proposed to avoid a costly initialization procedure to find a primal feasible starting point. Based on a standalone C~code implementation, we illustrate the computational performance of \software{PRESAS} against current state of the art QP solvers for multiple linear and nonlinear MPC case studies. We also show that the solver is real-time feasible on a dSPACE MicroAutoBox-II rapid prototyping unit for vehicle control applications, and numerical reliability is illustrated based on experimental results from a testbench of small-scale autonomous vehicles.}
	
	
	\jnlcitation{\cname{%
			\author{Quirynen R.}, and
			\author{Di Cairano S.}} (\cyear{2019}), 
		\ctitle{Block-Structured Preconditioning of Iterative Solvers within an Active-Set Strategy for Model Predictive Control Applications}, \cjournal{Optimal Control Applications and Methods}, \cvol{2019;00:1--6}.}
	
	\maketitle
	
		


\section{Introduction}

Optimization based control and estimation techniques, such as model predictive control~(MPC) and moving horizon estimation~(MHE), have experienced an increasing amount of interest from industry~\cite{Rawlings2017}. One of the main advantages of these methods is their systematic way of incorporating system dynamics and inequality constraints. At each sampling instant, MPC solves a multi-stage dynamic optimization problem that minimizes a particular cost function subject to continuity conditions and inequality constraints. A block-sparse quadratic program~(QP) structure arises in linear or linear time-varying formulations of predictive control and estimation. A similarly structured quadratic program~(QP) forms also the subproblem within a sequential quadratic programming~(SQP) method for nonlinear optimal control~\cite{Diehl2005,Gros2016}.

In this paper, we are interested in solving the following convex constrained linear-quadratic optimal control problem~(OCP)
\begin{subequations}
	\label{eq:QP_OCP}
	\begin{alignat}{3}
	&\underset{X,U}{\min}  &  \sum_{k=0}^{N-1} &\frac{1}{2} \begin{bmatrix} x_k\\u_k \end{bmatrix}^\top \begin{bmatrix}
	Q_k & S_k^\top \\ S_k & R_k
	\end{bmatrix} \begin{bmatrix} x_k\\u_k \end{bmatrix} + \begin{bmatrix} q_k\\r_k \end{bmatrix}^\top \begin{bmatrix} x_k\\u_k \end{bmatrix} &&+ \frac{1}{2} x_N^\top Q_N x_N + q_N^\top x_N \label{eq:QP_OCP:obj} \\
	& \;\st \quad &x_0 &= \bx, \label{eq:QP_OCP:init} \\
	&&x_{k+1} &= a_k + A_k x_k + B_k u_k, && k=0,\ldots,N{-}1, \label{eq:QP_OCP:dyn} \\
	&&0 &\geq d_k + D_{k}^\x x_k + D_{k}^\u u_k, && k=0,\ldots,N,  \label{eq:QP_OCP:path}
	\end{alignat}
\end{subequations}
where we define the state vectors as $x_k \in \R^{\nx}$, the control inputs as $u_k \in \R^{\nU}$, the \change{positive semidefinite} cost matrices as $Q_k \in \R^{\nx \times \nx}$, $S_k \in \R^{\nU \times \nx}$ and $R_k \in \R^{\nU \times \nU}$, and $N$ is the horizon length. We denote the state and control trajectory, respectively, as $X := [x_0^\top, \ldots, x_N^\top]^\top$ and $U := [u_0^\top, \ldots, u_{N-1}^\top]^\top$. The constraints include the system dynamics with $A_k \in \R^{\nx \times \nx}$, $B_k \in \R^{\nx \times \nU}$, the inequality constraints with $D_{k}^\x \in \R^{\nck \times \nx}$, $D_{k}^\u \in \R^{\nck \times \nU}$ and $\bx \in \R^{\nx}$ denotes the current state estimate.

A real-time implementation of MPC requires the solution of the block-structured QP in~\eqref{eq:QP_OCP} within a specified time period, typically on embedded control hardware with limited computational resources and a relatively small amount of available memory~\cite{DiCairano2018tutorial}. These challenges have sparked considerable efforts in research on embedded optimization algorithms that are tailored to predictive control applications such as described, e.g., in~\cite{Ferreau2017}.
The following requirements are important to be taken into account when designing or choosing an embedded QP solver for industrial applications of predictive control:
\begin{enumerate}
	\itemsep0em 
	\item scaling of computational complexity and memory requirements with problem dimensions $N$, $\nx$ and $\nU$,
	\item deterministic or early termination of solver in real-time applications~(obtaining a feasible but suboptimal solution),
	\item computational performance on embedded control hardware with relatively limited available resources, 
	\item warm starting capabilities to reduce solver execution time in receding horizon control applications,
	\item numerical reliability and ease of understanding, implementing and maintaining by non-experts,
	\item portability of solver code and limited software dependencies for embedded hardware.
\end{enumerate}
\change{Note that the latter list could additionally include infeasibility detection as a requirement for embedded MPC implementations. Tailored infeasibility detection for our proposed \software{PRESAS} QP solver is part of future research and remains outside the scope of the present paper. Instead, in what follows, we consider a soft constraint reformulation of all state-dependent inequality constraints, such that the real-time solver can provide a control solution at each time step. Note that input constraints always need to be respected as hard constraints instead. Based on the solution from the embedded optimization algorithm, another component in the overall control system could perform fault detection and recovery.
}

Most implementations of embedded optimization algorithms that have been successfully applied for real-time optimal control rely on direct linear algebra, e.g., in \software{qpOASES}~\cite{Ferreau2014}, \software{CVXOPT}~\cite{Andersen2013}, \software{FORCES}~\cite{Domahidi2014} or \software{HPMPC}~\cite{Frison2014}. However, it is known that iterative methods can result in a better asymptotic complexity when solving the saddle point linear systems arising in second order optimization methods~\cite{Benzi2005}. Iterative solvers, such as the minimal residual~(MINRES) or the conjugate gradient~(CG) method~\cite{Greenbaum1997}, are suitable for hardware acceleration of the linear system solution, e.g., using an FPGA for fast embedded applications~\cite{Boland2008}, due to their higher ratio of addition and multiplication operations. The use of an iterative method also allows to trade off between computational time and solution accuracy~\cite{Knyazev2015}. For a general linear system, iterative solvers tend to converge rather poorly without preconditioning~\cite{Benzi2005}.

Earlier work in~\cite{Shahzad2010b,Shahzad2010a} studied the use of block-diagonal preconditioners, in combination with numerical techniques that are tailored to an inexact interior-point~(IP) framework. \change{More recently, the work in~\cite{Otta2019} presents an active-set strategy based on a sparse variant of the Newton projection with proportioning~(NPP) algorithm and using an augmented Lagrangian-based preconditioner for the MINRES method.}
Unlike IP methods, an active-set quadratic programming algorithm can considerably benefit from the use of warm or hot-starting techniques to reduce the computational effort when solving a sequence of closely related optimal control problems, see~\cite{Bartlett2000,Ferreau2017,Wright1996}. 
\change{The cost per iteration of an active-set QP~solver is computationally cheaper} by exploiting low-rank updates of the matrix factorizations when changing the current guess for the active set~\cite{Ferreau2014,Kirches2011}, \change{compared to standard IP methods~\cite{Bartlett2000} or compared to the NPP algorithm in~\cite{Otta2019}}. We propose block-structured preconditioning techniques within a primal active-set strategy for real-time optimal control, following our initial investigation in~\cite{Quirynen2018a,Quirynen2018c}. The resulting Preconditioned RESidual solver within an Active-Set optimization strategy~(\software{PRESAS}) allows for an initial setup computational complexity of $\O(N m^3)$ and a per-iteration complexity of $\O(N m^2)$ for the QP in~\eqref{eq:QP_OCP}, where $m$ denotes the number of state and control variables in the system.

In the present paper, we provide a complete overview of our recent developments on block-structured preconditioning of iterative solvers with a primal active-set strategy, including block-diagonal preconditioners for MINRES as well as constraint preconditioning within a projected CG method. We describe theoretical properties for each of the preconditioning techniques, which are validated based on numerical simulation case studies in both single- and double-precision arithmetics. In addition, with respect to our preliminary studies in~\cite{Quirynen2018a,Quirynen2018c}, in this manuscript we also introduce an improved warm starting of the \software{PRESAS} solver, to avoid a costly procedure to find a primal feasible starting point, based on an augmented Lagrangian technique. Using an efficient and self-contained C~code implementation, the computational performance of \software{PRESAS} is illustrated against state of the art QP solvers. We show that the \software{PRESAS} QP~solver is real-time feasible on a dSPACE MicroAutoBox-II rapid prototyping unit for vehicle control applications, and numerical reliability is illustrated based on experimental results from a testbench of small-scale autonomous vehicles.

The paper is organized as follows. Section~\ref{sec:prelim} presents the preliminaries on embedded quadratic programming, primal active-set methods and on preconditioning of iterative solvers. Section~\ref{sec:prec} describes two block-diagonal preconditioners that can be used within a preconditioned minimal residual~(PMINRES) solver. Then, Section~\ref{sec:PRESAS_ppcg} presents and discusses an optimal control structured constraint preconditioner for the projected preconditioned conjugate gradient~(PPCG) method. Two warm-started initialization procedures for the \software{PRESAS} solver are presented in Section~\ref{sec:PRESAS}, and the resulting software implementation for embedded MPC applications is described in Section~\ref{sec:software}. The numerical properties and computational performance of different variants of the proposed \software{PRESAS} solver are illustrated in Section~\ref{sec:simResults} based on both linear and nonlinear MPC case studies, including hardware-in-the-loop simulations for vehicle control on a dSPACE MicroAutoBox-II. Section~\ref{sec:hamsters} presents experimental results for \software{PRESAS} running on a testbench of small-scale autonomous vehicles and Section~\ref{sec:concl} concludes the paper.

\section{Preliminaries} \label{sec:prelim}

We assume that the convex QP in~\eqref{eq:QP_OCP} has a unique global solution that is non-degenerate.
A solution is considered to be degenerate when either the strict complementarity condition or the linear independence constraint qualification~(LICQ) does not hold~\cite{Nocedal2006}.
In order to have a unique solution to the QP in~\eqref{eq:QP_OCP}, the Hessian needs to be positive definite on the null-space of the active constraints at the optimal solution. \change{For example, if $S_k = \zero$, the weighting matrices are chosen to satisfy $Q_k \succeq 0$ and $R_k \succ 0$ in common MPC formulations, such that the reduced Hessian is automatically positive definite for Eq.~\eqref{eq:QP_OCP}.}

\subsection{Embedded QP Solvers for Optimal Control}
\label{sec:convex}

There is a general trade-off between solvers that make use of second-order information and require only few but computationally expensive iterations, e.g., \software{qpOASES}~\cite{Ferreau2014}, versus first-order methods that are of low complexity but may require many more iterations, such as \software{PQP}~\cite{Cairano2013}, \software{ADMM}~\cite{Raghunathan2015} and other gradient or splitting-based methods~\cite{Ferreau2017}. In addition, there is an important distinction between optimal control algorithms that target the dense versus the sparse problem formulation. The numerical elimination of the state variables in a condensing routine~\cite{Bock1984,Frison2013a} is typically of a computational complexity $\O(N^2 m^3)$ but, unlike nonlinear problems, it can be mostly avoided in linear MPC applications. However, even with an offline preparation of the dense QP formulation, solvers applied to this dense QP will have a runtime complexity of $\O(N^2 m^2)$~\cite{Kirches2011}. \change{Examples of QP solvers to treat the condensed MPC formulation include \software{PQP}~\cite{Cairano2013}, \software{ADMM}~\cite{Raghunathan2015} and \software{qpOASES}~\cite{Ferreau2014}. Instead, we focus on directly solving the non-condensed MPC formulation with the block sparsity structure in~\eqref{eq:QP_OCP}, similar to the sparse solvers in \software{FORCES}~\cite{Domahidi2014}, \software{HPMPC}~\cite{Frison2014} and \software{qpDUNES}~\cite{Frasch2015}.}

It is important to note that many tailored QP algorithms for real-time optimal control rely on strict convexity of the cost function. This enables using a \change{sparsity exploiting} dual Newton strategy such as in \software{qpDUNES}~\cite{Frasch2015}, linear algebra routines such as the block-tridiagonal Cholesky factorization of the Schur complement in~\cite{Anderson1999,Wang2010a} or the particular Riccati recursion for linear-quadratic control problems in~\cite{Frison2013,Wright1996}. In the case of a positive semidefinite cost matrix, regularization needs to be applied, followed by an iterative refinement procedure to obtain a solution to the original problem~\cite{Frison2013}. 
This combination of regularization and iterative refinement can be also needed in the presence of ill-conditioned QP matrices. Instead, here, we do not assume strict convexity of the cost function and we propose to use an iterative method to solve each linear KKT system. 

\subsection{Primal Feasible Active-Set Method}

The basic idea behind active-set methods is to find an optimal active set by iteratively updating a current guess. When fixing the active inequality constraints at the current solution guess, a corresponding structured equality constrained QP is solved to compute a new search direction
\begin{subequations}
	\label{eq:EQP_OCP}
	\begin{alignat}{3}
	&\underset{\Delta X,\Delta U}{\min}  &  \sum_{k=0}^{N-1} &\frac{1}{2} \Delta w_k^\top H_k\, \Delta w_k + \tilde{h}_k^\top \Delta w_k &&+ \frac{1}{2} \Delta x_N^\top Q_N \Delta x_N + \tilde{q}_N^\top \Delta x_N \label{eq:EQP_OCP:obj} \\
	& \;\;\st \; &\Delta x_0 &= 0, \label{eq:EQP_OCP:init} \\
	&&\Delta x_{k+1} &= A_k \Delta x_k + B_k \Delta u_k, && k=0,\ldots,N{-}1, \label{eq:EQP_OCP:dyn} \\
	&&0 &= D_{k,i}^\x \Delta x_k + D_{k,i}^\u \Delta u_k,  &&({k,i}) \in \W,  \label{eq:EQP_OCP:path}
	\end{alignat}
\end{subequations}
where $\W$ denotes the current guess for the active set, i.e., the \emph{working set}. The variables $\Delta w_{k} := (\Delta x_{k},\Delta u_{k}) = (x_{k} - \bar{x}_{k},u_{k} - \bar{u}_{k})$ are defined for $k=0,\ldots,N-1$ and $\Delta w_N := \Delta x_N = x_N - \bar{x}_N$, \change{such that $\Delta X := [\Delta x_0^\top, \ldots, \Delta x_N^\top]^\top$ and $\Delta U := [\Delta u_0^\top, \ldots, \Delta u_{N-1}^\top]^\top$,} where $\bar{w}_{k} := (\bar{x}_{k},\bar{u}_{k})$ denotes the current guess for the optimal solution of the QP in~\eqref{eq:QP_OCP}. \change{Let us define the Lagrangian for the constrained QP in~\eqref{eq:QP_OCP} as follows
\begin{equation}
\begin{aligned}
\LL(X, U, \lambda_0, \mu_0, \ldots, \lambda_N, \mu_N) := &\sum_{k=0}^{N-1} \frac{1}{2} \begin{bmatrix} x_k\\u_k \end{bmatrix}^\top \begin{bmatrix}
Q_k & S_k^\top \\ S_k & R_k
\end{bmatrix} \begin{bmatrix} x_k\\u_k \end{bmatrix} + \begin{bmatrix} q_k\\r_k \end{bmatrix}^\top \begin{bmatrix} x_k\\u_k \end{bmatrix} + \frac{1}{2} x_N^\top Q_N x_N + q_N^\top x_N + \lambda_0^\top \left(\bx - x_0\right) \\
&+ \sum_{k=0}^{N-1} \lambda_{k+1}^\top \left(a_k + A_k x_k + B_k u_k - x_{k+1}\right) + \sum_{k=0}^{N} \mu_k^\top (d_k + D_{k}^\x x_k + D_{k}^\u u_k),
\end{aligned} \label{eq:Lagrangian}
\end{equation}
where $\lambda_k$ and $\mu_k$ denote the Lagrange multipliers, respectively, for equality and inequality constraints at each stage $k=0, \ldots, N$.} The Hessian matrices $H_k := \begin{bmatrix}
Q_k & S_k^\top \\ S_k & R_k
\end{bmatrix}$ and Lagrangian gradients $\tilde{h}_k := H_k \bar{w}_{k} + \begin{bmatrix} q_k - \lambda_{k}\\r_k \end{bmatrix} + \begin{bmatrix} A_k^\top\\B_k^\top \end{bmatrix} \lambda_{k+1} + \sum_{i \in \W_k}\mu_{k,i}\begin{bmatrix} D_{k,i}^{\x^\top}\\D_{k,i}^{\u^\top} \end{bmatrix}$ are defined for $k=0,\ldots,N-1$ and $\tilde{q}_N := Q_N \bar{x}_N + q_N - \lambda_{N} + \sum_{i \in \W_N}\mu_{N,i} D_{N,i}^{\x^\top}$.

The equality constrained QP in~\eqref{eq:EQP_OCP} results in the search direction $\bar{w}_{k} \,+\, \alpha \Delta w_{k}$ for which all constraints in the set $\W$ remain satisfied, regardless of the value $\alpha \ge 0$.
A distinction should be made between primal, dual and primal-dual active-set methods~\cite{Nocedal2006}. In addition, parametric methods have been proposed~\cite{Ferreau2014} in order to exploit the parametric aspect within a homotopy framework.
The active-set identification suffers from a combinatorial complexity and is therefore \emph{NP hard} in the worst case~\cite{Bartlett2000}. However, active-set methods have been successfully used in many real-time control applications~\cite{Ferreau2017}. 
In this work, we consider the general framework of primal feasible active-set methods as described in~\cite{Nocedal2006}, see Algorithm~\ref{alg:prim_AS}, and which generally provides the following numerical properties and advantages over alternative active-set optimization strategies:
\begin{enumerate}
	\itemsep0em 
	\item does not require strong duality and is relatively easy to implement compared to, e.g., a primal-dual method,
	\item can be terminated at any time point and still provides a feasible, yet suboptimal, solution to the OCP,
	\item automatically maintains the linear independence of the constraint gradients in its working set~\cite{Nocedal2006}.
\end{enumerate}
The latter property is very important in practice since it avoids the need for computing the maximal linearly independent subset of the current active set as described in~\cite{Ferreau2014} for dense QP problems. In case of sparse optimal control problems, the matrix factorization update techniques for active-set methods become relatively complicated, as discussed in~\cite{Kirches2011}. Also, it becomes nontrivial or impossible to efficiently maintain the linear independence of the working set in combination with a direct exploitation of the block sparsity structure in optimal control problems. Therefore, we focus on a primal feasible active-set method where such linear independence is automatically preserved.
\change{
The work in~\cite{Axehill2019,Axehill2020} describes an approach for complexity certification to compute a partitioning of the parameter space and corresponding bounds on the number of iterations in primal active-set methods, which could be applied directly to our proposed optimization algorithm for MPC.
Procedures for handling degeneracy and cycling in primal active-set methods can be found in~\cite{Nocedal2006}, but such issues are not common in practical MPC problems. 
Note that the primal active-set method in Algorithm~\ref{alg:prim_AS} adds~(step~$6$) or removes~(step~$11$) exactly one constraint in the current working set at each iteration. Multiple constraints with negative Lagrange multiplier values could instead be removed in step~$11$, but we restrict to the standard implementation in Algorithm~\ref{alg:prim_AS} for the remainder of this article.
}


\subsection{Block-Structured Linear System Solution}

At each iteration of the active-set method, the solution of the equality constrained QP in~\eqref{eq:EQP_OCP} corresponds to solving a saddle point linear system that represents the first order necessary conditions of optimality, i.e., the Karush-Kuhn-Tucker~(KKT) conditions.
The block-structured saddle point linear system can be written in compact form as follows
\begin{equation}
\begin{bmatrix} \H & \A^\top \\ \A & \zero \end{bmatrix} \bvec \Delta y \\ \Delta \nu \evec = \bvec b \\ 0 \evec \quad \text{or} \quad \K \,z = \bar{b}, \label{eq:linSys}
\end{equation}
where the primal variables $\Delta y := [\Delta w_0^\top, \ldots, \Delta w_N^\top]^\top$ and the Lagrange multipliers $\Delta \nu := [\Delta \lambda_0^\top, \Delta \mu_0^\top, \ldots, \Delta \lambda_N^\top, \Delta \mu_N^\top]^\top$ are defined. Similarly, the right-hand side vector is defined as $b := \change{-} [\tilde{h}_0^\top, \ldots, \tilde{h}_{N-1}^\top, \tilde{q}_N^\top]^\top$ based on the Lagrangian gradient vectors.
The linear system in~\eqref{eq:linSys} has a particular sparsity structure because $\H$ is the block-diagonal Hessian matrix and the constraint matrix $\A$ reads as
\begin{equation}
\begin{aligned}
\H &= \begin{bmatrix}
H_0 &  \\
 & H_1 & \\
 & & \ddots &  \\
&&& Q_N
\end{bmatrix}, \qquad
\A &= \begin{bmatrix}
-\eye & \zero \\
E_0^\x & E_0^\u & \\
A_0 & B_0 & - \eye & \zero \\
&&& \ddots
\end{bmatrix}
= \begin{bmatrix}
\begin{bmatrix} -\eye & \zero \end{bmatrix} \\
E_0 & \\
C_0 & \begin{bmatrix} -\eye & \zero \end{bmatrix} \\
& \ddots
\end{bmatrix},
\end{aligned}
\label{eq:jacMat}
\end{equation}
where \change{$\eye$ and $\zero$ denote, respectively, the identity and zero matrix,} $E_{k}^\x := [D_{k,i}^\x]_{i \in \W_k}$ and $E_{k}^\u := [D_{k,i}^\u]_{i \in \W_k}$ denote the active inequality constraints for each interval $k$, corresponding to the working set in $\W$, and we define the block matrices $C_k := \begin{bmatrix} A_k & B_k \end{bmatrix}$ and $E_k := \begin{bmatrix} E_{k}^\x & E_{k}^\u \end{bmatrix}$ for the sake of simplifying the notation.

\begin{algorithm}[h]%
	\caption{Primal Feasible Active-Set Method for Quadratic Programming.}
	\label{alg:prim_AS}
	\begin{algorithmic}[1]
		\Require Primal feasible starting point $y_0$ and working set $\W_0$.
		\For{$\text{iter}=1,2,\ldots,\text{max}$}
		\State \parbox[t]{\dimexpr\linewidth-\algorithmicindent}{Solve the QP in~\eqref{eq:EQP_OCP}, i.e., the linear system in~\eqref{eq:linSys}, to compute the search direction $(\Delta y,\Delta \nu)$.\strut}
		\State \change{$\alpha \gets \min_{({k,j}) \notin \W} -\frac{d_{k,j} + D_{k,j} \bar{w}_k}{D_{k,j} \Delta w_k}$. \hfill(step length computation)}
		\If{\change{$\alpha < 1$}}
		\State $(k,j) \gets \argmin_{\change{({k,j}) \notin \W}} -\frac{d_{k,j} + D_{k,j} \bar{w}_k}{D_{k,j} \Delta w_k}$. \hfill\change{(most blocking constraint)}
		\State $\W \gets \W \cup \{({k,j})\}$. \hfill(add constraint to working set)
		\State $y^\star \gets y + \alpha\,\Delta y$ and $\nu^\star \gets \nu + \alpha\,\Delta \nu$.
		\Else
		\State $y^\star \gets y + \Delta y$ and $\nu^\star \gets \nu + \Delta \nu$.
		\If{$\exists ({k,j}) \in \W: \mu_{k,j}^\star < 0$}
		\State $\W \gets \W \setminus \{ \argmin_{({k,j}) \in \W} \;\mu_{k,j}^\star \}$. \hfill(remove from working set)
		\Else
		\State {\bf stop} \hfill(QP solution found)
		\EndIf
		\EndIf
		\EndFor
		\Ensure Optimal solution $y^\star$ and $\nu^\star$ for the QP in~\eqref{eq:QP_OCP}.
	\end{algorithmic}
\end{algorithm}

\subsection{Preconditioning of Iterative Solvers}

In Eq.~\eqref{eq:linSys}, the matrix $\A$ has full rank and $\H$ is symmetric and positive semidefinite. Unlike prior work on embedded optimization algorithms for optimal control based on direct linear algebra routines, e.g.,~\cite{Frison2013,Kirches2011,Wang2010a,Wright1996}, we propose the use of iterative solvers as discussed for general saddle point linear systems in~\cite{Benzi2005,Benzi2008}. However, preconditioning is necessary for the good performance of iterative solvers~\cite{Knyazev2015}. \change{Standard preconditioning results in a modified linear system $\T^{-1} \K \,z = \T^{-1} b$ where $\T$ is the preconditioner. Instead, based on a positive definite preconditioning matrix $\T = L L^\top$, symmetric preconditioning $L^{-1} \K\,L^{-\top}$ can be performed, in which the preconditioner $\T$ should be designed such that}
\begin{enumerate}
	\itemsep0em 
	\item computations with the operator \change{$\T = L L^\top$} are cheaper than solving the original saddle point linear system in~\eqref{eq:linSys},
	\item the preconditioned matrix \change{$L^{-1} \K\,L^{-\top}$} approximates the identity or its eigenvalues are tightly clustered~\cite{Greenbaum1997}.	
\end{enumerate}

An overview on algebraic and application-specific preconditioners can be found in~\cite{Benzi2005}. Here, we focus on two block-diagonal preconditioning techniques and one constraint preconditioner that is tailored to the projected conjugate gradient method.
We aim to utilize warm- or hot-starting capabilities of active-set methods and exploit the rank-one property of active-set changes in the updating procedure for the factorization of the structured preconditioner from one active-set iteration to the next. 



\section{Optimal Control Structured Block-Diagonal Preconditioning}
\label{sec:prec}

Next, we briefly present two positive definite block-diagonal preconditioners for solving each optimal control structured linear system in~\eqref{eq:linSys} by using the minimal residual~(MINRES) method. Both preconditioning techniques have been studied extensively for solving general saddle point linear systems, e.g., in~\cite{Benzi2008}.


\subsection{Augmented Lagrangian~(AL) Preconditioner}
\label{sec:PRESAS-AL}
To accelerate convergence for the iterative solver, we can use the block-diagonal preconditioner
\begin{equation}
\Ta = \begin{bmatrix} \H + \A^\top \V \A & \zero \\ \zero & \V^{-1} \end{bmatrix} = \begin{bmatrix} \H + \gamma \A^\top \A & \zero \\ \zero & {\gamma}^{-1}\, \eye \end{bmatrix},  \label{eq:block_pred_Aug}
\end{equation}
where $\V$ is a symmetric positive definite weighting matrix. A popular choice for the weighting matrix, which follows an augmented Lagrangian type approach~\cite{Benzi2008}, is $\V = \gamma\, \eye$ where $\gamma > 0$ is a scalar and $\eye$ denotes the identity matrix. The application of the augmented Lagrangian type preconditioner $\Ta$ in~\eqref{eq:block_pred_Aug} for optimal control requires the factorization of the block-tridiagonal sparse matrix $\H + \gamma \A^\top \A$,
\begin{equation}
\begin{aligned}
&\H + \gamma \A^\top \A = \begin{bmatrix}
\hat{H}_0 + \gamma\, G_0^\top G_0 & - \gamma\, C_0^\top \\
- \gamma\, C_0 & \hat{H}_1 + \gamma\, G_1^\top G_1 & - \gamma\, C_1^\top \\
& - \gamma\, C_1 &  & \ddots \\
&& \ddots  & \ddots & - \gamma\, C_{N-1}^\top \\
&&& - \gamma\, C_{N-1} & \hat{H}_N + \gamma\, E_N^\top E_N
\end{bmatrix},
\end{aligned} \label{eq:augMat}
\end{equation}
where $G_k = \begin{bmatrix} E_k \\ C_k \end{bmatrix} \change{= \begin{bmatrix} E_{k}^\x & E_{k}^\u \\ A_k & B_k \end{bmatrix}}$ \change{and} $\hat{H}_k = \begin{bmatrix} Q_k + \gamma\, \eye & S_k^\top \\ S_k & R_k \end{bmatrix}$.

The structured matrix in Eq.~\eqref{eq:augMat} is positive definite for any value $\gamma > 0$ so that the block-tridiagonal Cholesky decomposition~\cite{Anderson1999} can be applied. The asymptotic computational complexity for this matrix factorization is $\O(N m^3)$, where $m = \nx+\nU$. However, such computational cost is incurred only once per QP solution for the initial guess of active inequality constraints. At each iteration of the active-set method, in which one constraint is added or removed from the working set, a corresponding row is added or removed from the Jacobian matrix $\A$ in~\eqref{eq:jacMat}. Given this modification of the matrix $\H + \gamma \A^\top \A$, a rank-one update to its block-tridiagonal Cholesky factorization can be computed at a computational complexity of $\O(N m^2)$ in each iteration of the solver. In this work, we further rely on a reverse or backward implementation of such a block-tridiagonal Cholesky decomposition, as proposed also in~\cite{Frasch2015} for the \software{qpDUNES} solver. Under the assumption that active-set changes are often more likely to occur near the beginning of the prediction horizon, rank-one factorization updates can typically be performed more efficiently for such a reverse or backward Cholesky decomposition.

\subsection{Schur Complement~(SC) Preconditioner}
\label{sec:PRESAS-SC}
Another approach is based on using a Schur complement type preconditioner
\begin{equation}
\Ts = \begin{bmatrix} \Hp & \zero \\ \zero & \A \Hp^{-1} \A^\top \end{bmatrix},  \label{eq:block_pred_Schur}
\end{equation}
where $\Hp \approx \H$ such that $\Hp \succ 0$. \change{For example, the approximate Hessian matrix $\Hp$ can be defined based on a simple diagonal regularization $\Hp = \H + \epsilon\, \eye \succ 0$.} The preconditioner $\Ts$ in~\eqref{eq:block_pred_Schur} results in a block-tridiagonal sparsity structure 
\begin{equation}
\begin{aligned}
\A \Hp^{-1} \A^\top 
	&= \change{\begin{bmatrix}
	\Ix \tilde{H}_0^{-1} \Ix^\top & -\Ix \tilde{H}_0^{-1} G_0^\top &  \\
	-G_0 \tilde{H}_0^{-1} \Ix^\top & G_0 \tilde{H}_0^{-1} G_0^\top + \begin{bmatrix} \zero & \zero \\ \zero & \Ix \tilde{H}_1^{-1} \Ix^\top \end{bmatrix} & -\begin{bmatrix} \zero \\ \Ix \tilde{H}_1^{-1} G_1^\top \end{bmatrix} &   \\
	& -\begin{bmatrix} \zero & G_1 \tilde{H}_1^{-1}\Ix^\top \end{bmatrix} & G_1 \tilde{H}_1^{-1} G_1^\top + \begin{bmatrix} \zero & \zero \\ \zero & \Ix \tilde{H}_2^{-1} \Ix^\top \end{bmatrix} & \ddots \\
	&&\ddots&\ddots & -\tilde{Q}_N^{-1} E_N^{\x^\top} \\
	&&& -E_N^\x \tilde{Q}_N^{-1} & E_N^\x \tilde{Q}_N^{-1} E_N^{\x^\top}
	\end{bmatrix},}
\end{aligned} \label{eq:schurMat}
\end{equation}
\change{where the matrix $\Ix := \begin{bmatrix} \eye & \zero \end{bmatrix}$ is defined to simplify the notation.}
Each of the block matrices in the preconditioning matrix $\Ts$ in Eq.~\eqref{eq:schurMat} is of different dimensions, corresponding to the number of active inequality constraints in each block. The asymptotic computational complexity for the block-tridiagonal Cholesky factorization of the matrix $\A \Hp^{-1} \A^\top$, as part of the preconditioner in~\eqref{eq:block_pred_Schur}, is $\O(N m^3)$ in which $m = \nx+\nU$.
Similar to before, a rank-one factorization update can be computed, resulting in a computational complexity of $\O(N m^2)$ at each iteration of the active-set method.

In case that a particular Hessian block matrix $H_k$ is positive semidefinite, a positive definite approximation $\tilde{H}_k = \begin{bmatrix} \tilde{Q}_k & \tilde{S}_k^\top \\ \tilde{S}_k & \tilde{R}_k \end{bmatrix} \succ 0$ needs to be computed. One can apply an \emph{on-the-fly} regularization in the form of a modified Cholesky factorization, e.g., as in~\cite{Frasch2015,Nocedal2006}, or alternatively the structure-exploiting regularization technique tailored to optimal control that was recently proposed in~\cite{Verschueren2017}. Here, we consider a simple but standard regularization procedure of the form $\tilde{H}_k = H_k + \epsilon\, \eye$, where the value for $\epsilon > 0$ is chosen sufficiently small such that $\tilde{H}_k \approx H_k$ forms a good Hessian approximation, but large enough such that $\tilde{H}_k \succ 0$ is positive definite and leads to good numerical conditioning of $\Ts$ in~\eqref{eq:schurMat}. 

\subsection{Preconditioned MINRES Algorithm}

Both the Schur-complement based $\Ts$ and the augmented Lagrangian preconditioner $\Ta$ are symmetric positive definite, \change{given any positive value for $\epsilon > 0$ or $\gamma > 0$}. 
Therefore, these preconditioning techniques can be used within the preconditioned minimal residual~(PMINRES) algorithm as described in~\cite{Greenbaum1997}. PMINRES is based on the three-term recurrence in the Lanczos iteration for symmetric matrices and specifically requires a positive definite preconditioner~\cite{Benzi2005}. On the other hand, the generalized minimal residual~(GMRES) algorithm is based on the Arnoldi iteration, for which the computational cost generally grows with each iteration, but it does not need a positive definite preconditioner. 
Here, we focus on using PMINRES in combination with either of the block-diagonal preconditioners in Eq.~\eqref{eq:block_pred_Aug} or~\eqref{eq:block_pred_Schur}.

As discussed in detail by~\cite{Greif2006}, as the value of $\gamma > 0$ increases for the augmented Lagrangian preconditioner \change{$\Ta = \La \La^\top$}, the eigenvalues of the preconditioned matrix \change{$\La^{-1} \K\,\La^{-\top}$} become more and more tightly clustered around $\pm 1$. For a sufficiently large value of $\gamma$, MINRES can therefore converge within two iterations in the ideal setting, but choosing $\gamma$ too large may result in ill-conditioning of the matrix $\Ta$. 
For the Schur-complement based preconditioner \change{$\Ts = \Ls \Ls^\top$}, when $\Hp=\H$ is invertible, the preconditioned matrix has three distinct eigenvalues $1$, $\frac{1}{2}(1+\sqrt{5})$ and $\frac{1}{2}(1-\sqrt{5})$, and therefore MINRES converges within three iterations~\cite{Murphy2000}. The eigenvalues become clustered when, e.g., regularization is needed to make $\Hp \approx \H$ positive definite~\cite{Benzi2005,Murphy2000}. Similar to before, for a smaller value of $\epsilon > 0$ in $\Hp = \H + \epsilon\, \eye$, the eigenvalues of \change{$\Ls^{-1} \K\,\Ls^{-\top}$} become more tightly clustered. Yet, when choosing the value of $\epsilon > 0$ too small, ill-conditioning of the matrix $\Ts$ may cause numerical issues. The effect of these choices on the overall performance of the algorithm will be illustrated in the case studies of Section~\ref{sec:simResults}.

%

\section{Projected Preconditioned Conjugate Gradient Method}
\label{sec:PRESAS_ppcg}

The parameters, $\gamma$ or $\epsilon$, in the above mentioned preconditioners can be difficult to choose, especially when using lower precision arithmetics or in the presence of ill-conditioned QP matrices. Therefore, in what follows, we present a preconditioned iterative solver for which the performance is less dependent on a careful choice for the regularization parameter. We propose to use conjugate gradient~(CG) iterations in the null space of the active inequality constraints, based on a particular projection operator that preserves the block-structured sparsity in optimal control.


\subsection{Conjugate Gradient for Reduced Linear System}

The equality constrained QP in~\eqref{eq:EQP_OCP} could be solved by eliminating the constraints and solving the resulting reduced problem formulation. 
Next, we briefly review the issues that arise when applying the CG method to such a reduced formulation of the linear KKT system. We define the matrix $\Z$ as a basis for the null space of the constraint matrix $\A$ such that $\A \Z = 0$. 
Given the linear system in~\eqref{eq:linSys}, the reduced formulation then reads as
\begin{equation}
\Z^\top \H \Z \,\Delta y_z = \Z^\top b, \label{eq:red_linSys}
\end{equation}
where $\Delta y = \Z \,\Delta y_z$, i.e., the solution vector $\Delta y$ is in the null space of the active constraint matrix $\A$. 

One could directly apply the conjugate gradient~(CG) method with standard preconditioning techniques to the reduced linear system in~\eqref{eq:red_linSys} as described in~\cite{Coleman2001,Gould2001a}. This is generally not advisable for optimal control problems because the block-structured sparsity is destroyed, similar to the use of a condensing routine~\cite{Bock1984}, resulting in expensive CG iterations of complexity $\O(N^2 m^2)$ instead of the desired runtime complexity of $\O(N m^2)$. In addition, forming a null space basis matrix $\Z$, computing the reduced Hessian $\Z^\top \H \Z$ and constructing and applying a preconditioner for this reduced Hessian is computationally expensive~\cite{Gould2001a}, in general. Sparsity exploiting update procedures for the null space basis matrix and for the factorization in a Schur complement step of an active-set strategy have been proposed in~\cite{Kirches2011}, but they are relatively complex to implement.

\subsection{Projected Preconditioned Conjugate Gradient~(PPCG)}

Instead of reducing to the linear system in Eq.~\eqref{eq:red_linSys}, we propose to apply a projection operator $\P$, onto the null space of $\A$, to the linear KKT system in Eq.~\eqref{eq:linSys}:
\begin{equation}
\P\, \H \,\Delta y = \P\, b, \label{eq:proj_linSys}
\end{equation}
where $\P \A^\top = 0$~\cite{Bakhvalov1994}. One possible definition of the projection operator is $\P := \Z \left(\Z^\top \Z\right)^{-1} \Z^\top$, but it requires again the computation of a basis matrix $\Z$ for the null space of $\A$. The idea behind the projected preconditioned conjugate gradient~(PPCG) method from~\cite{Coleman2001,Dollar2005,Gould2001a} is to use the alternative projection operator
\begin{equation}
\P_{\A:\Hp} := \eye - \A^\top\left(\A \Hp^{-1} \A^\top\right)^{-1} \A \Hp^{-1}, \label{eq:proj_A}
\end{equation}
where $\eye$ denotes the identity matrix. The projection in~\eqref{eq:proj_A} does not require any basis matrix $\Z$ and preserves the sparsity of a block-structured problem, \change{due to the block-tridiagonal sparsity structure of the positive definite matrix $\A \Hp^{-1} \A^\top$ and the block-diagonal structure of the positive definite matrix $\Hp$}, as discussed in the next subsection. The projection operator in~\eqref{eq:proj_A} is defined for a particular preconditioner $\Hp \succ 0$ and $\Hp \approx \H$ that can be applied directly to Eq.~\eqref{eq:proj_linSys}, resulting in
\begin{equation}
\Hp^{-1} \, \P_{\A:\Hp} \, \H \,\Delta y = \Hp^{-1} \, \P_{\A:\Hp} \, b. \label{eq:pred_proj_linSys}
\end{equation}
One can then apply the conjugate gradient method to the projected preconditioned linear system in~\eqref{eq:pred_proj_linSys}.

The operation $\left(\Hp^{-1} \, \P_{\A:\Hp}\right) b_1$ in Eq.~\eqref{eq:pred_proj_linSys} can alternatively be computed by solving the linear system
\begin{equation}
\begin{bmatrix} \Hp & \A^\top \\ \A & \zero \end{bmatrix} \bvec z_1 \\ z_2 \evec = \bvec b_1 \\ 0 \evec, \quad \text{or} \quad \Tc \,z = \bar{b}_1, \label{eq:constrPrec}
\end{equation}
for a certain vector $b_1$. The matrix $\Tc$ is also referred to as a \emph{constraint preconditioner} for the KKT system in~\eqref{eq:linSys}. The solution of the linear system~\eqref{eq:constrPrec} corresponds to
\begin{equation}
z_1 := \Hp^{-1} \left( b_1 - \A^\top z_2 \right), \; z_2 := \left(\A \Hp^{-1} \A^\top\right)^{-1} \A \Hp^{-1} b_1, \label{eq:constrPrec2}
\end{equation}
such that $z_1 = \left(\Hp^{-1} \, \P_{\A:\Hp}\right) b_1$ holds. Furthermore, we use the improved variant of the PPCG method based on the residual update strategy in~\cite{Dollar2005,Gould2001a}, as described in Algorithm~\ref{alg:ppcg}, in order to avoid significant round-off errors in the range of $\A^\top$. An additional iterative refinement step could also be used to improve the numerical stability of the solver.

\begin{algorithm}[h]
	\caption{Projected Preconditioned Conjugate Gradient~(PPCG) method with the residual update strategy.}
	\label{alg:ppcg}
	\begin{algorithmic}[1]
		\Require Matrices $\H$, $\A$, Hessian approximation $\Hp$ and vector $b$ in~\eqref{eq:linSys}.
		\State $\Delta y \gets 0$ and $r \gets \H\, \Delta y - b$.
		\State Solve $\begin{bmatrix} \Hp & \A^\top \\ \A & \zero \end{bmatrix} \bvec g \\ v \evec = \bvec r \\ 0 \evec$ such as in Eq.~\eqref{eq:constrPrec2}.
		\State $p \gets -g$, $\Delta \nu \gets -v$ and $r \gets r - \A^\top v$.
		\For {$n=1,\ldots,n_{\mathrm{max}}$}
		\State $\alpha \gets \frac{r^\top g}{p^\top H\, p}$.
		\State $\Delta y \gets \Delta y + \alpha p$. 
		\State $r^+ \gets r + \alpha \H\, p$. 
		\If {$\Vert r^+\Vert < \epsilon_{\mathrm{tol}}$}
		\State $\Delta y^* \gets \Delta y$ and $\Delta \nu^* \gets \Delta \nu$.
		\State {\bf stop} \hfill(solution found)
		\Else
		\State Solve $\begin{bmatrix} \Hp & \A^\top \\ \A & \zero \end{bmatrix} \bvec g^+ \\ v^+ \evec = \bvec r^+ \\ 0 \evec$ such as in Eq.~\eqref{eq:constrPrec2}.
		\State $\Delta \nu \gets \Delta \nu - v^+$.
		\State $\beta \gets \frac{r^{+^\top} g^+}{r^\top g}$.
		\State $p \gets -g^+ + \beta p$.
		\State $g \gets g^+$.
		\State $r \gets r^+ - \A^\top v^+$.
		\EndIf
		\EndFor
		\Ensure Solution vectors $\Delta y^*$ and $\Delta \nu^*$ for linear system in Eq.~\eqref{eq:linSys}.
	\end{algorithmic}
\end{algorithm}

\subsection{Block-structured Sparse Constraint Preconditioner}

The standard block-tridiagonal sparsity structure of the positive definite matrix $\A \Hp^{-1} \A^\top$ as well as the block-diagonal structure of the positive definite matrix $\Hp$ can be exploited in Eq.~\eqref{eq:constrPrec2} for the implementation of \software{PRESAS} in combination with the PPCG method for fast MPC applications. In addition, the approximate Hessian matrix $\Hp$ in the preconditioner is often chosen to be diagonal, in practice. 
As described earlier in Section~\ref{sec:prec}, the block-tridiagonal Cholesky factorization for the matrix $\A \Hp^{-1} \A^\top$ can be efficiently maintained based on a rank-one factorization update for each active-set change. Similar to the techniques in Section~\ref{sec:prec}, this variant of the solver requires an initial setup computational complexity of $\O(N m^3)$ and a per-iteration complexity of $\O(N m^2)$ for solving the optimal control structured QP in~\eqref{eq:QP_OCP}.

The quality of the constraint preconditioner is defined by the accuracy of the Hessian approximation $\Hp \approx \H$ such that $\Hp \succ 0$. However, unlike the case for the Schur complement type block-diagonal preconditioner, the eigenvalues of the preconditioned matrix
can be shown to be
\begin{equation}
\{ 1 \} \;\cup\; \sigma\left( \left(\Z^\top \Hp \Z\right)^{-1} \Z^\top \H \Z\right), \label{eq:spectrum}
\end{equation}
where $\sigma(\cdot)$ denotes the spectrum of a matrix. As described in~\cite{Gould2001a}, typical choices are $\Hp = \eye$ or $\Hp = \texttt{diag}(\H)$. Here, we use again a regularized Hessian approximation of the form $\Hp = \H + \epsilon\, \eye$, where $\epsilon > 0$ is chosen sufficiently small without causing ill-conditioning of the constraint preconditioner.

Two conclusions can be made straightforwardly from the expression for the eigenvalues of the preconditioned matrix in Eq.~\eqref{eq:spectrum}. First, all eigenvalues are equal to $1$ in case that $\H = \Hp \succ 0$, i.e., only one CG iteration is needed. Second, all eigenvalues remain equal to $1$ also when an augmented Lagrangian type regularization is applied to the Hessian, i.e.,  $\Hp = \H + \A^\top \V \A \succ 0$ such that the reduced Hessian reads $\Z^\top \Hp \Z = \Z^\top \H \Z$ as in~\cite{Dollar2005}. The latter is very common for MPC applications, e.g., when performing a slack reformulation of inequality constraints based on an exact L1 penalty. The full Hessian matrix will be positive semidefinite as a result of such slack reformulation, but a regularization of the slack variables in the constraint preconditioner does not change the reduced Hessian and therefore it does not affect the convergence of the PPCG method.
In practice, the number of iterations can be larger than one for both of these cases because of numerical round-off errors that are caused by ill-conditioning and/or the usage of low-precision arithmetics.

\section{Warm-started Initialization of \software{PRESAS} for Predictive Control}
\label{sec:PRESAS}
	
As indicated in Algorithm~\ref{alg:prim_AS}, a primal active-set method requires an initial point that is primal feasible. In the general case of a constrained quadratic program, this is a nontrivial task that corresponds to a \emph{Phase~I} procedure as described, for example, in~\cite{Fletcher1987,Nocedal2006}. Here, we describe two approaches that are tailored to the solution of constrained optimal control problems, in order to allow an infeasible start of the primal active-set method and therefore avoid the need for an explicit \emph{Phase~I} procedure.

\subsection{Forward Simulation of Dynamics and Slack Reformulation}
\label{sec:PRESAS_sim}

When solving the parametric OCP in~\eqref{eq:QP_OCP} within receding horizon based control or estimation, it becomes relatively easy to satisfy the initial value condition and the continuity constraints based on a forward simulation of the state variables using a shifted version of the control trajectory. In addition, one may introduce slack variables in order to always satisfy all state-dependent inequality constraints in~\eqref{eq:QP_OCP:path}. 
The first approach is therefore based on a penalty method, similar in spirit to that described in~\cite{Nocedal2006}. Let us reformulate the structured QP obtained from optimal control problem~\eqref{eq:QP_OCP} as follows:
\begin{subequations}
	\label{eq:soft_QP_OCP}
	\begin{alignat}{3}
	&\underset{X,U,\S}{\min}  &  \sum_{k=0}^{N-1} &\frac{1}{2} \begin{bmatrix} x_k\\u_k\\\s_k \end{bmatrix}^\top \begin{bmatrix}
	Q_k & S_k^\top \\ S_k & R_k \\ && M_k
	\end{bmatrix} \begin{bmatrix} x_k\\u_k\\\s_k \end{bmatrix} + \begin{bmatrix} q_k\\r_k\\m_k \end{bmatrix}^\top \begin{bmatrix} x_k\\u_k\\\s_k \end{bmatrix} &&+ \frac{1}{2} \begin{bmatrix} x_N\\\s_N \end{bmatrix}^\top \begin{bmatrix}
	Q_N &  \\  & M_N
	\end{bmatrix} \begin{bmatrix} x_N\\\s_N \end{bmatrix} + \begin{bmatrix} q_N\\m_N \end{bmatrix}^\top \begin{bmatrix} x_N\\\s_N \end{bmatrix} \label{eq:soft_QP_OCP:obj} \\
	& \;\st \quad &x_0 &= \bx, \label{eq:soft_QP_OCP:init} \\
	&&x_{k+1} &= a_k + A_k x_k + B_k u_k, && k=0,\ldots,N{-}1, \label{eq:soft_QP_OCP:dyn} \\
	&&\s_{k} &\geq d_k + D_{k}^\x x_k + D_{k}^\u u_k, && k=0,\ldots,N,  \label{eq:soft_QP_OCP:path} \\
	&&\s_{k} &\geq 0, && k=0,\ldots,N,  \label{eq:soft_QP_OCP:slack}
	\end{alignat}
\end{subequations}
where $M_k \ge 0$ and $m_k \ge 0$, respectively, represent the L2 and exact L1 penalty for each of the slack variables $\S := [\s_0, \ldots, \s_{N}]^\top$. It can be shown that $\s_{k} = 0$ for each $k=0,\ldots,N$ in the solution of the soft-constrained QP~\eqref{eq:soft_QP_OCP}, whenever a feasible solution exists for the QP in~\eqref{eq:QP_OCP} and for a sufficiently large value of $m_k > 0$ in the exact L1 penalty~\cite{Fletcher1987}.

The resulting solver is easy to implement, because there is no explicit \emph{Phase~I} procedure needed for initialization and a primal active-set method automatically maintains the linear independence of the constraint gradients in its working set~\cite{Nocedal2006}. However, the forward simulation procedure to satisfy the equality constraints in eqs.~\eqref{eq:soft_QP_OCP:init}-\eqref{eq:soft_QP_OCP:dyn} may lead to a poor initial guess for the primal active-set method in case of unstable system dynamics and/or large changes of the initial state value in~\eqref{eq:soft_QP_OCP:init}, after shifting the state and control trajectory from one control time step to the next.

\subsection{Augmented Lagrangian Method for Initial Value Condition}
\label{sec:PRESAS_WSA}

In order to avoid the potential numerical instability in the forward simulation procedure due to unstable system dynamics and/or large changes of the initial state value, we propose to enforce the initial state value condition of Eq.~\eqref{eq:soft_QP_OCP:init} in the form of an augmented Lagrangian method as described in~\cite{Nocedal2006}. For the soft-constrained optimal control structured problem in~\eqref{eq:soft_QP_OCP}, we remove the initial state value constraint and instead compute the solution to the following QP
\begin{subequations}
	\label{eq:alm_soft_QP_OCP}
	\begin{alignat}{3}
	&\underset{X,U,\S}{\min}  &  \sum_{k=0}^{N-1} &\frac{1}{2} \begin{bmatrix} x_k\\u_k\\\s_k \end{bmatrix}^\top \begin{bmatrix}
	Q_k & S_k^\top \\ S_k & R_k \\ && M_k
	\end{bmatrix} \begin{bmatrix} x_k\\u_k\\\s_k \end{bmatrix} + \begin{bmatrix} q_k\\r_k\\m_k \end{bmatrix}^\top \begin{bmatrix} x_k\\u_k\\\s_k \end{bmatrix} &&+ \frac{1}{2} \begin{bmatrix} x_N\\\s_N \end{bmatrix}^\top \begin{bmatrix}
	Q_N &  \\  & M_N
	\end{bmatrix} \begin{bmatrix} x_N\\\s_N \end{bmatrix} + \begin{bmatrix} q_N\\m_N \end{bmatrix}^\top \begin{bmatrix} x_N\\\s_N \end{bmatrix} \label{eq:alm_soft_QP_OCP:obj1} \\
	&&&+ \frac{1}{2} \Vert x_0 - \bx \Vert_{\Rho^i}^2 + \lambda_0^{i^\top} \left( \bx - x_0 \right) \label{eq:alm_soft_QP_OCP:obj2} \\
	& \;\st \quad &x_{k+1} &= a_k + A_k x_k + B_k u_k, && k=0,\ldots,N{-}1, \label{eq:alm_soft_QP_OCP:dyn} \\
	&&\s_{k} &\geq d_k + D_{k}^\x x_k + D_{k}^\u u_k, && k=0,\ldots,N,  \label{eq:alm_soft_QP_OCP:path} \\
	&&\s_{k} &\geq 0, && k=0,\ldots,N,  \label{eq:alm_soft_QP_OCP:slack}
	\end{alignat}
\end{subequations}
at each iteration $i = 0, \ldots$, where $\lambda_0^{i}$ and $\Rho^i := \texttt{diag}(\rho^i)$, respectively, denote the Lagrange multipliers and the quadratic penalty for the initial value condition of Eq.~\eqref{eq:soft_QP_OCP:init}. As discussed in~\cite{Nocedal2006}, after solving the QP~\eqref{eq:alm_soft_QP_OCP}, the Lagrange multiplier values for the initial state value condition can be updated by
\begin{equation}
\lambda_0^{i+1} = \lambda_0^{i} + \Rho^i \,\left( \bx - x_0 \right).\label{eq:mult_update}
\end{equation}
An increasingly large value can be chosen for the quadratic penalty $\Rho^{i+1} \succeq \Rho^i$, even though convergence of the augmented Lagrangian method can be assured even without increasing the value for $\Rho := \texttt{diag}(\rho)$~\cite{Nocedal2006}. In addition, each quadratic penalty can be updated individually for each of the initial state vector components.

\begin{algorithm}[h]
	\caption{Warm-starting of \software{PRESAS}: Augmented Lagrangian Method for Initial State Value~(\software{PRESAS-WSA}).}
	\label{alg:aug_lag}
	\begin{algorithmic}[1]
		\Require Initial solution guess $X^0$, $U^0$, $\S^0$ and initial values $\Rho^0 = \texttt{diag}(\rho^0)$ and $\lambda^0$.
		\State Forward simulation of system dynamics to ensure satisfaction of Eq.~\eqref{eq:soft_QP_OCP:dyn}.
		\State Increase slack variables $\s_k \ge 0$ for $k=0, \ldots, N$ to satisfy inequality constraints~\eqref{eq:soft_QP_OCP:path}.
		\While {$\Vert x_0 - \bx \Vert > \epsilon_{\mathrm{tol}}$}
		\State Solve block-structured QP subproblem in~\eqref{eq:alm_soft_QP_OCP} using Algorithm~\ref{alg:PRESAS}~(\software{PRESAS}), with reuse of Cholesky factorization.
		\State Update Lagrange multiplier values: $\lambda_0^{i+1} \gets \lambda_0^{i} + \Rho^i \,\left( \bx - x_0 \right)$.
		\For {$j=1,\ldots,\nx$}
			\If {$| x_{0,j} - \hat{x}_{0,j} | > \epsilon_{\mathrm{tol}}$ AND $\rho_j^i < \rho_{\mathrm{max}}$}
			\State Update quadratic penalty value: $\rho_j^{i+1} \ge \rho_j^i$ for initial state variable $x_{0,j}$.
			\If {$\rho_j^{i+1} > \rho_j^i$}
			\State Rank-1 Cholesky factorization update for preconditioner in \software{PRESAS}.
			\EndIf
			\EndIf
		\EndFor
		\EndWhile
		\Ensure Solution variables $X^\star$, $U^\star$ and $\S^\star$ for the soft-constrained QP in~\eqref{eq:soft_QP_OCP}.
	\end{algorithmic}
\end{algorithm}

Algorithm~\ref{alg:aug_lag} provides a detailed description of the augmented Lagrangian type implementation for effective warm starting of the proposed primal active-set strategy in the \software{PRESAS} QP solver, which will be referred to as \software{PRESAS-WSA}. \change{Note that the QP solution in each iteration is used as the warm-started initialization of \software{PRESAS} in the next augmented Lagrangian type iteration. In addition, the block-structured Cholesky factorization for the preconditioning matrix can be reused from one iteration to the next, resulting in a \emph{hot start} for each call to the \software{PRESAS} QP solver within the augmented Lagrangian method of Algorithm~\ref{alg:aug_lag}.} Rank-one factorization updates are needed for each of the individual changes in the quadratic penalty values, corresponding to a rank-one update for the Hessian matrix of the QP subproblem~\eqref{eq:alm_soft_QP_OCP}. However, as mentioned also earlier, when using a reverse block-tridiagonal Cholesky factorization procedure~\cite{Frasch2015}, the computational cost per rank-one update in the first stage of the control horizon becomes $\O(m^2)$ instead of the computational complexity of $\O(N m^2)$ for the complete factorization update. The computational performance of \software{PRESAS}, using either the penalty or the augmented Lagrangian method, is illustrated based on closed-loop numerical simulations for case studies of linear and nonlinear MPC in Section~\ref{sec:simResults}.

\section{Software Implementation of \software{PRESAS} for Embedded MPC Applications}
\label{sec:software}


The proposed \software{PRESAS} solver is described in Algorithm~\ref{alg:PRESAS}, based on the primal active-set strategy in Algorithm~\ref{alg:prim_AS} and the \software{PPCG} method in Algorithm~\ref{alg:ppcg}, in combination with structure exploiting factorization updates for the block-sparse preconditioner. The solver adheres to all embedded optimal control requirements~($1$-$6$) that were stated in the introduction. Namely, \software{PRESAS} enjoys the preferred computational complexity of $\O(N m^2)$ per iteration based on its block sparsity structure exploitation and the use of iterative linear algebra routines. For being an active-set method with tailored optimal control type structure exploitation, the proposed approach is considerably easier to implement than prior work such as~\cite{Kirches2011}. 

\begin{algorithm}[h]%
	\caption{Preconditioned RESidual method within primal Active-Set~(\software{PRESAS}) solver for MPC.}
	\label{alg:PRESAS}
	\begin{algorithmic}[1]
		\Require Primal feasible starting point $y_0$ and working set $\W_0$ (see warm starting in Section~\ref{sec:PRESAS}).
		\State Compute regularized Hessian block matrix $\tilde{H}_k \approx H_k$ such that $\tilde{H}_k \succ 0$ for $k=0,\ldots,N$.
		\If{$\tilde{H}_k$ is diagonal}
		\State Compute $\tilde{H}_k^{-1}$ by inverting diagonal elements for each $k=0,\ldots,N$. \hfill $\O(N m)$
		\Else
		\State Compute Cholesky factorization of block matrix $\tilde{H}_k$ for $k=0,\ldots,N$. \hfill $\O(N m^3)$
		\EndIf
		\State Block-tridiagonal (reverse) Cholesky factorization for matrix $\A \Hp^{-1} \A^\top$ in Eq.~\eqref{eq:constrPrec2}. \hfill $\O(N m^3)$
		\For{$\text{iter}=1,2,\ldots,\text{max}$}
		\State Algorithm~\ref{alg:ppcg} to solve the linear system in~\eqref{eq:linSys}, to compute the search direction $(\Delta y,\Delta \nu)$. \hfill $\O(N m^2)$
		\State \change{$\alpha \gets \min_{({k,j}) \notin \W} -\frac{d_{k,j} + D_{k,j} \bar{w}_k}{D_{k,j} \Delta w_k}$. \hfill(step length computation)}
		\If{\change{$\alpha < 1$}}
		\State $(k,j) \gets \argmin_{\change{({k,j}) \notin \W}} -\frac{d_{k,j} + D_{k,j} \bar{w}_k}{D_{k,j} \Delta w_k}$. \hfill\change{(most blocking constraint)}
		\State $\W \gets \W \cup \{({k,j})\}$. \hfill(add constraint to working set)
		\State Add row to matrix $\A$ and compute rank-one factorization update for $\A \Hp^{-1} \A^\top$. \hfill $\O(N m^2)$
		\State $y^\star \gets y + \alpha\,\Delta y$ and $\nu^\star \gets \nu + \alpha\,\Delta \nu$.
		\Else
		\State $y^\star \gets y + \Delta y$ and $\nu^\star \gets \nu + \Delta \nu$.
		\If{$\exists ({k,j}) \in \W: \mu_{k,j}^\star < 0$}
		\State $\W \gets \W \setminus \{ \argmin_{({k,j}) \in \W} \;\mu_{k,j}^\star \}$. \hfill(remove from working set)
		\State Remove row from $\A$ and compute rank-one factorization update for $\A \Hp^{-1} \A^\top$. \hfill $\O(N m^2)$
		\Else
		\State {\bf stop} \hfill(QP solution found)
		\EndIf
		\EndIf
		\EndFor
		\Ensure Optimal solution $y^\star$ and $\nu^\star$ for the QP in~\eqref{eq:QP_OCP}.
	\end{algorithmic}
\end{algorithm}


\subsection{Self-contained C~code Implementation for Embedded Applications}

In what follows, we illustrate the numerical properties and the computational performance of the proposed primal active-set solver for fast MPC applications, based on a preliminary C~code implementation. In order to simplify its use for industrial applications of optimization based control and estimation, the software implementation does not rely on any external libraries, e.g., for performing linear algebra routines. More specifically, all matrix factorizations and matrix-matrix multiplications are performed by standard triple loop implementations. It is well known that such textbook implementations can easily be outperformed, for medium- to large-scale matrix dimensions, by more advanced Basic Linear Algebra Subprogram~(BLAS) libraries. For this purpose, high-performance open-source BLAS implementations could be used instead such as \software{GotoBLAS}~\cite{Goto2008} and \software{OpenBLAS}~\cite{OpenBLAS2011}, as well as proprietary implementations such as Intel's \software{MKL} and AMD's \software{ACML}. The open-source \software{BLASFEO}~\cite{Frison2017b} framework provides hardware-tailored optimized dense BLAS routines, outperforming most alternative tools for small- to medium-scale matrix dimensions that occur typically for applications of embedded optimization. However, in what follows, we illustrate that our simplified self-contained C~code implementation remains very competitive with state of the art QP solvers. The main reasons for this are that \software{PRESAS} exploits the block-structured sparsity in direct optimal control problems and it does not require any matrix factorizations or matrix-matrix multiplications within the active-set iterations.

In addition, the C~code implementation is written such that the solver can easily be embedded in prototyping platforms such as, e.g., the dSPACE MicroAutoBox-II unit, as well as in microcontrollers that are used in real-world mechatronic systems. The computational hardware that we target for implementation has capabilities in the same scale of current automotive microcontrollers. As it is well known, due to the need to operate in harsh environmental conditions, the hard real-time requirement, and the cost considerations, current automotive microcontrollers are significantly less computationally powerful than standard desktop computers~\cite{DiCairano2018tutorial}, in terms of both memory and operations per second. This limits the amount of computations and data storage that we use in implementing our numerical optimization algorithms. Table~\ref{tab:PRESAS_variants} illustrates the variants of the proposed primal active-set method that are implemented in self-contained C~code. 
\change{Our implementation allows the flexibility of switching between different structure exploiting preconditioning techniques within the proposed block-sparse \software{PRESAS} solver, in order to adapt to a particular problem formulation.}
Table~\ref{tab:PRESAS_variants} also includes an implementation of our solver that is based on dense linear algebra, called \software{D-PRESAS}, which will be used in comparisons of the numerical simulation results in Section~\ref{sec:simResults}. By comparing \software{D-PRESAS} with \software{PRESAS}, based on the same preconditioned iterative solver, we can show the scalability of the tailored block-sparse structure exploitation for fast MPC implementations. As expected, \software{D-PRESAS} can remain competitive for relatively small-scale applications or it can even outperform the sparse \software{PRESAS} solver, and showing its capabilities allows us to assess the performance of the preconditioned iterative solver with respect to direct linear algebra methods.

\begin{table*}[tpb]
	\vspace*{-1pt}
	\caption{Most important variants of the primal active-set solver that have been implemented in self-contained C~code.}
	\vspace*{-6pt}
	\label{tab:PRESAS_variants}
	\centering
	\setlength{\tabcolsep}{0.3em}
	\begin{tabular}{l | c | c | c }
		\toprule
		\multicolumn{1}{c}{ }   	& \multicolumn{1}{c}{Sparsity structure}   	& \multicolumn{1}{c}{Preconditioned iterative solver}   	& \multicolumn{1}{c}{Primal-feasible warm starting}     	      \\
		\midrule
		\software{PRESAS-SC}     	& Block-sparse   	 &   PMINRES - Schur complement~(SC)   		&   Forward simulation + slack variables \\
		\software{PRESAS-AL}       	& Block-sparse   	 &   PMINRES - Augmented Lagrangian~(AL)   	&   Forward simulation + slack variables \\
		\software{PRESAS-PPCG}       	& Block-sparse   &   PPCG - Constraint preconditioner (Alg.~\ref{alg:ppcg})   		&   Forward simulation + slack variables \\
		\software{PRESAS-PPCG-WSA}      & Block-sparse   &   PPCG - Constraint preconditioner (Alg.~\ref{alg:ppcg})   		&   Augmented Lagrangian method (Alg.~\ref{alg:aug_lag}) \\
		\software{D-PRESAS-PPCG}       	& Dense  		 &   PPCG - Constraint preconditioner (Alg.~\ref{alg:ppcg})   		&   Slack variables \\
		\bottomrule
	\end{tabular}
	\vspace*{-6pt}
\end{table*}


%

\subsection{Nonlinear MPC: ACADO Code Generation Tool Interface}
\label{sec:NMPC_software}

\begin{figure}[tpb]
	\centerline{\hbox{
			\includegraphics[width=0.5\textwidth]{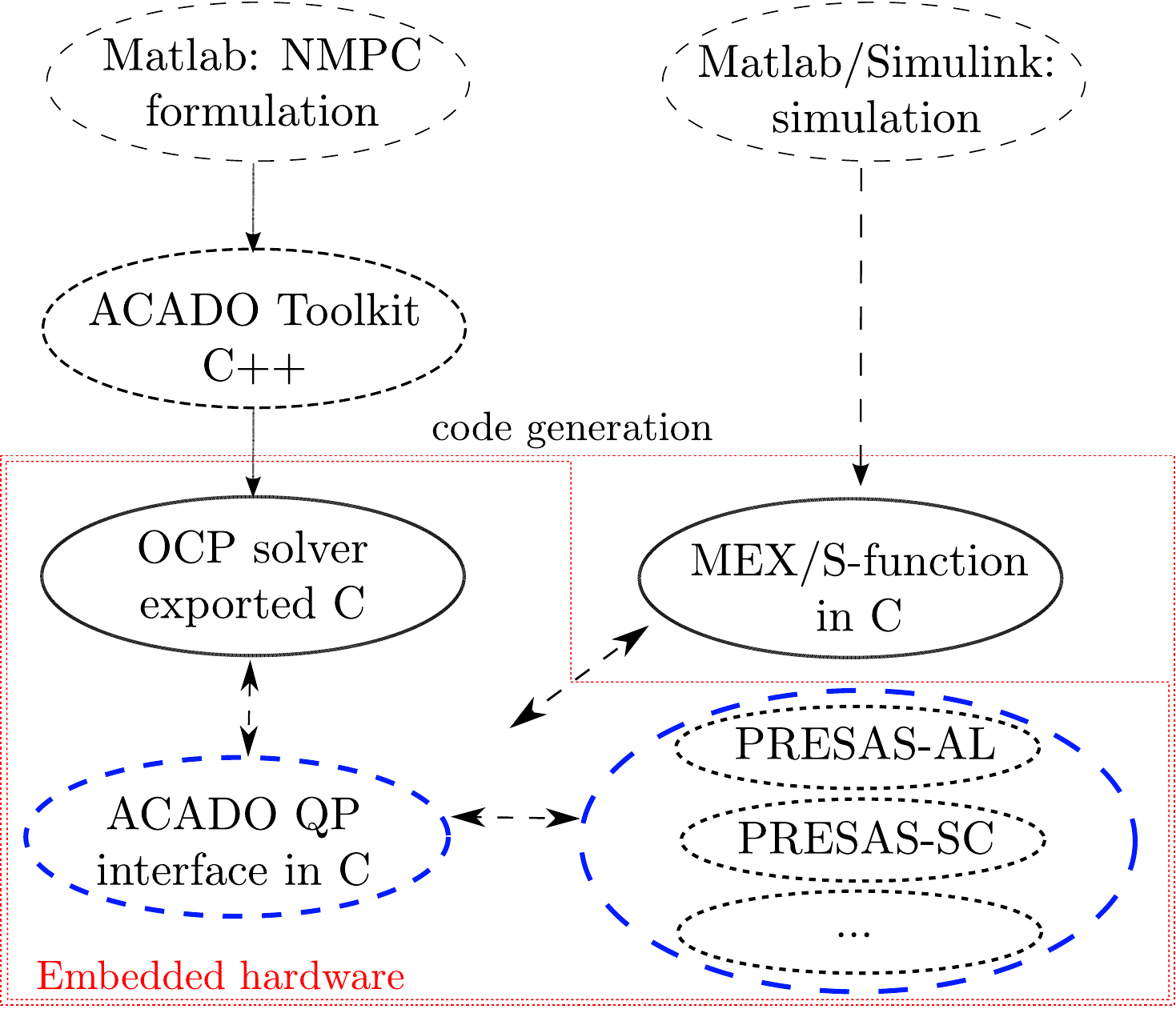}}}
	\caption{Illustration of the proposed work flow based on the open-source \software{ACADO Toolkit} software and its C~code generation tool for nonlinear MPC: custom sparse QP solver interface for \software{PRESAS}.}
	\label{fig:ACADO_PRESAS}
\end{figure}

In addition to linear and linear time-varying MPC applications, we illustrate the computational performance of the \software{PRESAS} solver also for nonlinear MPC~(NMPC) applications. For this purpose, we use a minimal software interface between the code generation tool of the \software{ACADO Toolkit}~\cite{Houska2011,Quirynen2014a} and the self-contained C~code implementations of \software{PRESAS}. The nonlinear optimal control solver in the open-source toolkit uses an online variant of a sequential quadratic programming~(SQP) type method, known as the real-time iterations~(RTI) scheme~\cite{Diehl2005}. The RTI approach is based on performing one SQP iteration per control time step, and on using a continuation-based warm starting of the state and control trajectories from one time step to the next. Each iteration consists of two steps:
\begin{enumerate}
	\itemsep0em 
	\item \emph{Preparation phase}: discretize and linearize the system dynamics, linearize the remaining constraint functions, and evaluate the quadratic objective approximation to build the optimal control structured QP subproblem of the form in~\eqref{eq:QP_OCP}.
	\item \emph{Feedback phase}: solve the local QP approximation~\eqref{eq:QP_OCP} to update the current values for all optimization variables, i.e., the state and control trajectories, and obtain the next control input to apply to the system.
\end{enumerate} 
Note that we developed two separate interfaces for either the block-sparse and dense versions of the \software{PRESAS} solver in Table~\ref{tab:PRESAS_variants}. The dense \software{D-PRESAS} variant can be used together with condensing~\cite{Frison2013a}, to numerically eliminate the state variables from the QP subproblem at each control time step, in the auto generated RTI optimization algorithm.
The resulting work flow, based on our custom \software{ACADO}-\software{PRESAS} software interface, is illustrated in Figure~\ref{fig:ACADO_PRESAS} and consists of the following three steps:
\begin{enumerate}
	\itemsep0em 
	\item \emph{Formulation}: high-level problem formulation in Matlab or C++ using ACADO syntax, including the continuous-time nonlinear system dynamics, least squares type objective function and inequality constraints.
	\item \emph{Generation}: automatic generation of C~code to implement an RTI algorithm for the user-specified optimal control problem. The generated code defines routines for a generalized Gauss-Newton~(GGN)~\cite{} method that performs tailored function evaluations and algorithmic differentiation~(AD). More specifically, the OCP solver performs numerical simulation and sensitivity propagation for discretization and linearization of the nonlinear system of differential equations, and it computes a linearization of the nonlinear inequality constraints and a Gauss-Newton type Hessian approximation.
	\item \emph{Compilation}: the generated \software{ACADO} code can be compiled together with the self-contained \software{PRESAS} solver and its custom interface, and linked into a MEX file or S-function such that it can be used from Matlab or Simulink.
\end{enumerate}


\section{Numerical Simulation Case Studies}
\label{sec:simResults}

We consider three different numerical simulation case studies. The first case study concerns an MPC problem formulation based on linear system dynamics for a chain of spring-connected masses as in~\cite{Wang2010a}. The second case study is a nonlinear dynamic system that describes an inverted pendulum on top of a cart. We consider the case of linear MPC for stabilization of the inverted pendulum and nonlinear MPC for the swing-up of the pendulum~\cite{Quirynen2014a}. Finally, a hardware-in-the-loop simulation of a trajectory tracking NMPC algorithm for vehicle control is presented in the third case study. We illustrate the computational performance for different variants of the \software{PRESAS} solver in Table~\ref{tab:PRESAS_variants}, based on standalone C~code implementations. 


\subsection{Case Study~1: MPC on Chain of Oscillating Masses}

This first test problem consists of the chain of oscillating masses, which is often used as a benchmark example for fast MPC solvers~\cite{Quirynen2018a,Wang2010a}. The linear time-invariant system dynamics and corresponding OCP formulation, of the form in~\eqref{eq:QP_OCP}, are described in more detail in~\cite{Wang2010a}. The full state of the system consists of the displacement and velocity of the $\nm$ masses, i.e., $x(t) \in \R^{2 \nm}$ such that the state dimension can be changed by changing the amount of masses. A number of actuators $\nU < \nm$ apply tensions between certain masses while respecting the actuator limitations as well as constraints on the position and velocity of each of the masses. In order to guarantee that the QP remains feasible at each sampling instant, the state constraints are softened by adding one slack variable for each control interval. Based on a sufficiently large penalization of this additional variable in the objective, a feasible solution can be found whenever possible. \change{As discussed in~\cite{Boyd2011}, the \software{ADMM} solver does not require the explicit use of slack variables in the OCP problem formulation to soften state-dependent inequality constraints based on soft thresholding. Therefore, unlike for all the other optimization algorithms, the numerical results in this paper are obtained without including auxiliary slack variables for \software{ADMM}.} During the closed-loop MPC simulations, a disturbance force, randomly generated with uniform distribution and kept the same for all simulations, acts on each of the masses as in~\cite{Wang2010a}.

\paragraph{Scaling Computational Complexity of \software{PRESAS} Solver with Problem Dimensions}

\begin{figure}[tpb]
	\centerline{\hbox{
			\includegraphics[width=0.95\textwidth]{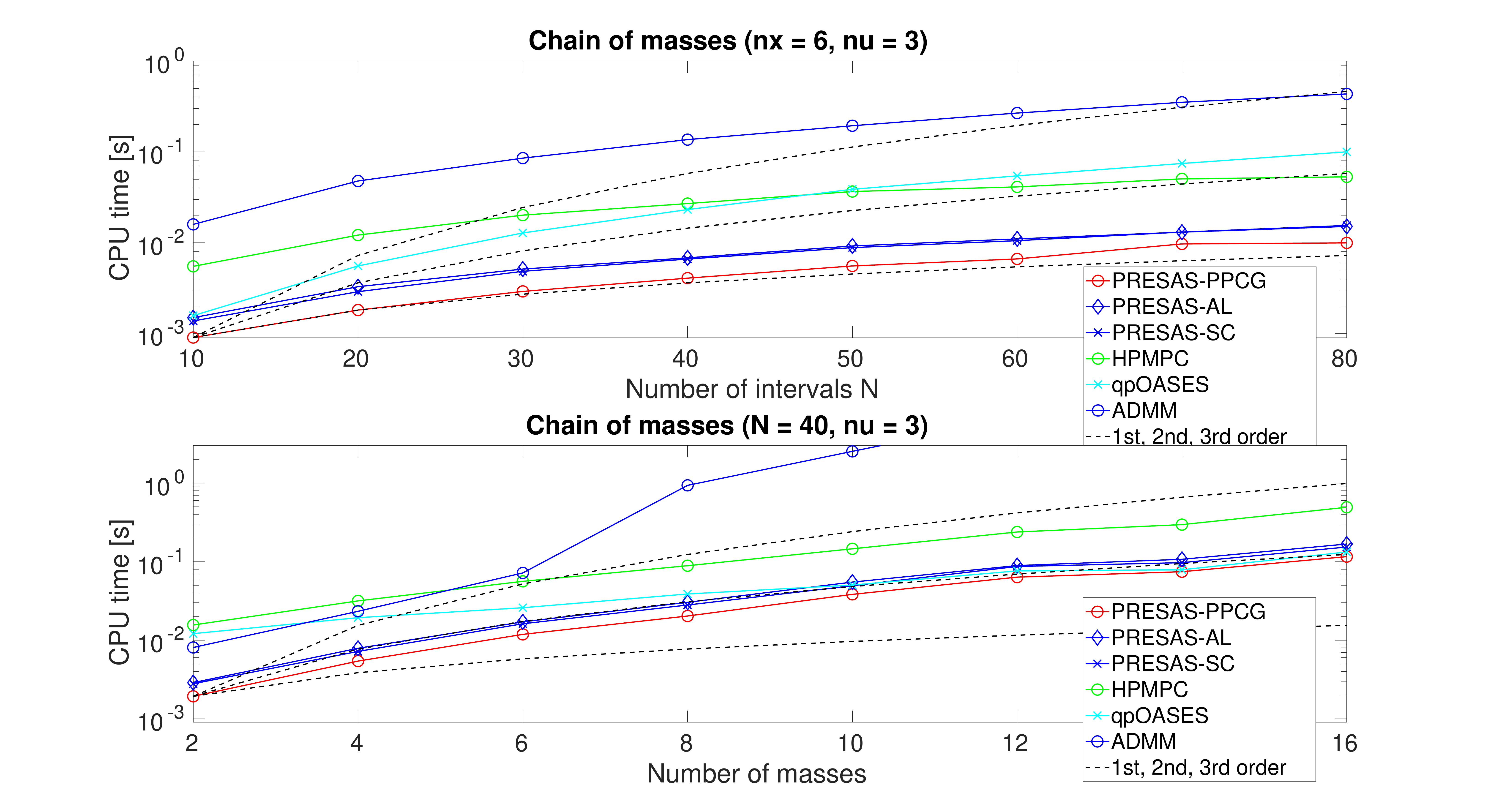}}}
	\caption{Average computation time~(ARM Cortex-A7 processor in the Raspberry Pi~2) of different QP solvers for linear MPC on the chain of masses with L2 slack reformulation: varying number of control intervals $N$ and number of masses $\nm$.}
	\label{fig:chain_results}
\end{figure}

Figure~\ref{fig:chain_results} shows the average computation time per QP solution during the closed-loop simulations of linear MPC for different numbers of masses $\nm$ and for a varying control horizon length $N$. \change{For simplicity, we do not consider the case where the control horizon is shorter than the prediction horizon, resulting in relatively less control variables, which could be used to reduce the computational cost of any optimization algorithm. Figure~\ref{fig:chain_results}} includes the computation time for \software{ADMM}, \software{qpOASES}, \software{HPMPC} and the $3$ variants of \software{PRESAS}. The considered \software{ADMM} algorithm~\cite{Raghunathan2015} and \software{qpOASES}~\cite{Ferreau2014} both solve the small but dense QP after numerically eliminating the state variables. 
The timing results for \software{ADMM} and \software{qpOASES} include the computation time for this condensing routine~\cite{Frison2013a}. Note that the black dashed lines in Figure~\ref{fig:chain_results} illustrate, respectively, a computational cost that increases with a $1^{\text{st}}$, $2^{\text{nd}}$ or $3^{\text{rd}}$ order of complexity, i.e., $\O(N)$, $\O(N^2)$ and $\O(N^3)$ in the top part, and $\O(m)$, $\O(m^2)$ and $\O(m^3)$ in the bottom part of Figure~\ref{fig:chain_results}. The typical runtime complexity of $\O(N^2 m^2)$ for the two dense solvers, \software{ADMM} and \software{qpOASES}, can be observed. However, note that the computation time for \software{qpOASES} remains nearly constant with respect to the number of masses, because the dense QP includes only the control variables and Figure~\ref{fig:chain_results} is constructed for a fixed amount of \change{$\nU=3$} control inputs. When increasing the number of masses in the bottom part of Figure~\ref{fig:chain_results}, both the ratio of control to state variables and the ratio of active inequality constraints decrease. On the other hand, the computational complexity of $\O(N m^3)$ for the interior-point method in \software{HPMPC}~\cite{Frison2014} and the per iteration complexity of $\O(N m^2)$ for the proposed variants of the \software{PRESAS} solver can be observed in Figure~\ref{fig:chain_results}.
Even though the \software{PRESAS-AL} and \software{PRESAS-SC} solvers are very competitive with state of the art optimal control algorithms, it can be observed that the \software{PRESAS-PPCG} method further outperforms these algorithms. This performance scales well with the number of control intervals and the number of state variables in this particular case study.

%

\paragraph{Solver Sensitivity to Parameter Selection and Numerical Accuracy}

\begin{figure}[tpb]
	\centerline{\hbox{
			\includegraphics[width=0.95\textwidth]{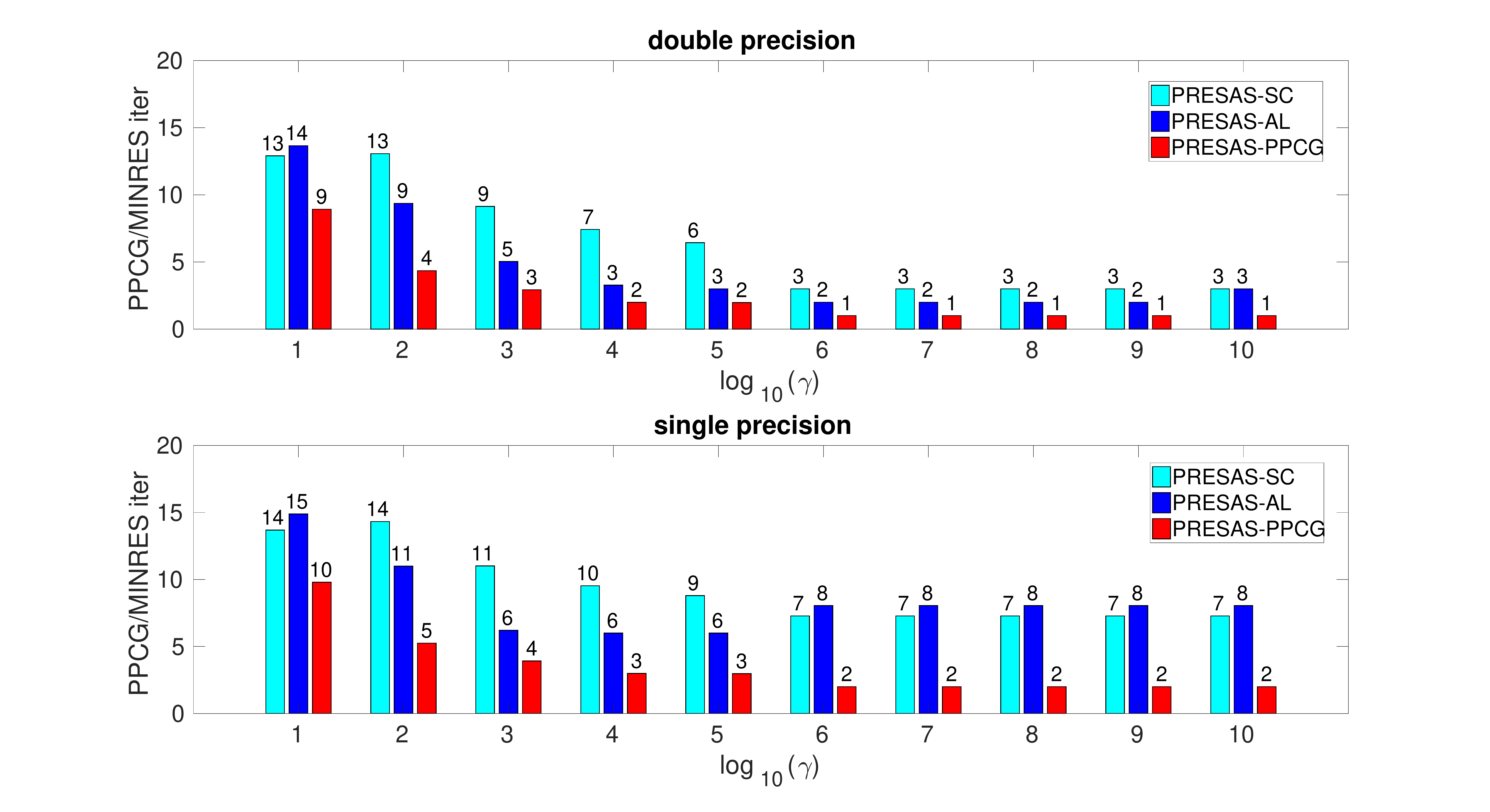}}}
	\caption{Average number of residual iterations for PPCG in \software{PRESAS-PPCG} or for PMINRES in \software{PRESAS-SC} or \software{PRESAS-AL} per linear system solution within the closed-loop MPC simulations for a chain of $\nm=3$ masses.}
	\label{fig:chain_results_gamma}
\end{figure}

\software{PRESAS-AL} requires a suitable choice for its parameter $\gamma$, while \software{PRESAS-SC} and the \software{PRESAS-PPCG} solver typically need diagonal regularization using a parameter $\epsilon = \frac{1}{\gamma}$. Figure~\ref{fig:chain_results_gamma} illustrates the numerical dependency of the number of PMINRES and PPCG iterations on the value of $\gamma$, either using single- or double-precision arithmetics. As expected, the value for $\gamma$ should be sufficiently large but not too large in order to maintain a good conditioning of the preconditioner.
In case of double-precision arithmetics, the parameter design appears to be much less sensitive. For a sufficiently large value of $\gamma$, the \software{SC}, \software{AL} and \software{PPCG} variant of \software{PRESAS}, respectively, require $3$, $2$ and $1$ iteration of the residual method on average for each solution of a saddle point linear system. When using single-precision arithmetics, the numerical performance is reduced for all three preconditioners. The \software{PRESAS-PPCG} method shows a consistently lower number of iterations \change{and it appears to be less sensitive to the choice of parameter value $\gamma$ for the considered range of values. In addition, the \software{PRESAS-PPCG} method} allows an average of only two \software{PPCG} iterations for a sufficiently small regularization parameter $\epsilon = \frac{1}{\gamma}$ in this case study. \change{In order to avoid ill conditioning of the preconditioner, a good choice for the parameter value is, e.g., $\gamma = 10^7$ for double precision and $\gamma = 10^4$ for single-precision arithmetics. Note that this parameter value is fixed for the numerical results that are presented in the remainder of this article.}

\subsection{Case Study~2: MPC for Inverted Pendulum on a Cart}

This second numerical case study consists of two closed-loop MPC simulations that, respectively, consider stabilization and swing-up of the inverted pendulum mounted on top of a cart~\cite{Quirynen2014a}. 
We present detailed timing results for an ARM Cortex-A7 processor in the Raspberry Pi~2. While they are not embedded processors by themselves, such Raspberry Pis use ARM cores with computational capabilities on the order of high-end microcontrollers that can be used for embedded real-time control applications in several industries~\cite{DiCairano2018tutorial}.

\paragraph{Warm-started Initialization for \software{PRESAS} solver: Inverted Pendulum Stabilization}

\begin{figure}[tpb]
	\centerline{\hbox{
			\includegraphics[width=0.95\textwidth]{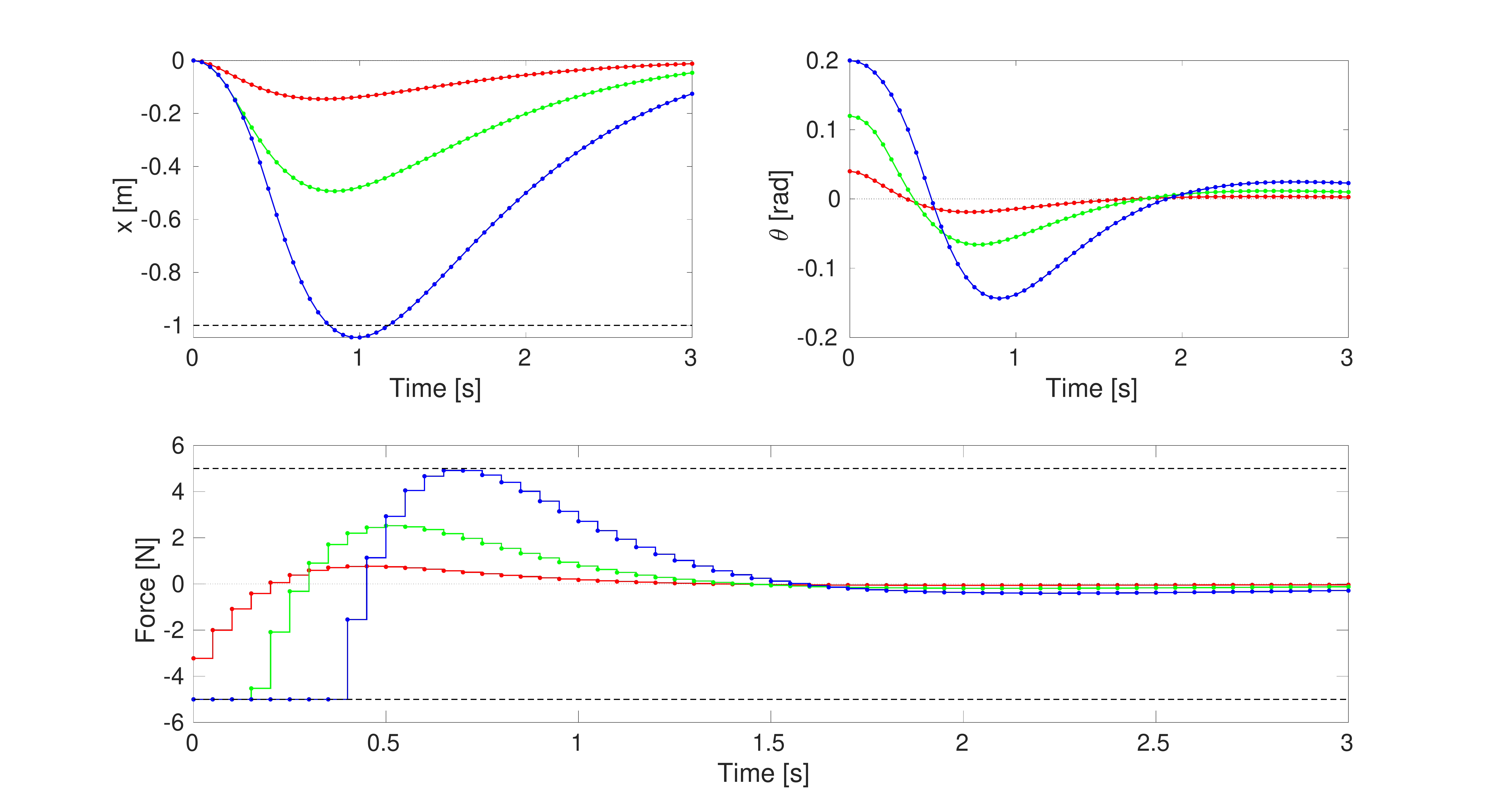}}}
	\vspace*{-8pt}
	\caption{Closed-loop trajectories for linear MPC to stabilize an inverted pendulum on top of a cart, starting from different initial conditions for the pendulum angle~($\theta_0=0.04$, $\theta_0=0.12$ or $\theta_0=0.20$~rad).}
	\vspace*{-8pt}
	\label{fig:pend_results}
\end{figure}

This first simulation involves stabilizing the nonlinear inverted pendulum in the upward unstable position, based on a linearization of the nonlinear dynamics in that steady state. The closed-loop trajectories are illustrated in Figure~\ref{fig:pend_results} for different initial conditions for the angle of the pendulum $\theta_0 = 0.04$,~$0.12$~or~$0.20$~rad. Both the actuated force and the cart position are constrained to remain within their respective bounds. As illustrated in Figure~\ref{fig:pend_results}, \change{it becomes impossible to satisfy the position constraint, resulting in a constraint violation around $t=1$~s in case of the initial value $\theta_0 = 0.20$, which is treated by a softened reformulation of the state-dependent inequality constraints.} This results in a relatively high number of online active-set changes with respect to the amount of state $\nx=4$ and control variables $\nU=2$. 

Table~\ref{tab:pend_results} shows the average and worst-case computation times on an ARM Cortex-A7 processor and the numbers of iterations for linear MPC of the inverted pendulum, including results for the QP solvers \software{ADMM}, \software{qpOASES}, \software{HPMPC}, \software{qpDUNES} and different variants of our proposed \software{PRESAS} solver.
For all the solvers, except for the interior point method in~\software{HPMPC}, warm starting is performed. For most of the QP solvers, the computational effort grows significantly when increasing the initial value for the angle of the inverted pendulum.
The results show that the \software{PRESAS-PPCG} algorithm consistently outperforms both the \software{PRESAS-SC} and the \software{PRESAS-AL} variants of our solver. In combination with the augmented Lagrangian based warm-started initialization procedure, see Section~\ref{sec:PRESAS_WSA}, the resulting \software{PRESAS-PPCG-WSA} solver demonstrates an improved worst-case computational performance. \change{Note that \software{qpOASES} reports $0$~iterations in Table~\ref{tab:pend_results} when the initial guess for the active set is optimal. And the results for the \software{PRESAS-PPCG-WSA} solver refer to the total number of inner iterations in Algorithm~\ref{alg:aug_lag}, i.e., it corresponds to the total number of active-set changes, which is equal to the total number of linear systems that are solved to compute each QP solution.} More specifically, for this particular case study, our solver is at least a factor $2$ times faster than all other state of the art QP solvers in Table~\ref{tab:pend_results}.
In addition, by using the primal active-set method in \software{PRESAS}, a feasible but suboptimal solution can always be obtained by limiting the maximum number of active-set changes. 
\begin{table*}[tpb]
\vspace*{-1pt}
\caption[]{Average and worst-case computation times~(ms) for MPC of an inverted pendulum~($\Ts = 50$~ms, $N=50$), with varying initial conditions for the angle of the pendulum, using double precision on an ARM Cortex-A7 processor in the Raspberry Pi~2.~\footnotemark}
\vspace*{-6pt}
\label{tab:pend_results}
\centering
\setlength{\tabcolsep}{0.9em}
\begin{tabular}{l | c c | c c | c c  }
	\toprule
	\multicolumn{1}{c}{ }   	& \multicolumn{2}{c}{$\theta_0=0.04$}   	& \multicolumn{2}{c}{$\theta_0=0.12$}   	& \multicolumn{2}{c}{$\theta_0=0.20$}     	      \\
	\midrule
	& time~[ms] & \# iter   					& time~[ms] & \# iter   					& time~[ms] & \# iter				 \\
	& (mean/max) & (mean/max)   				& (mean/max) & (mean/max)   				& (mean/max) & (mean/max)				 \\
	\midrule

\software{ADMM}  & $8.10/63.56$ & $17.9/92$   &   $17.24/90.51$ & $42.1/181$   &   $25.94/109.79$ & $65.4/239$ \\
\software{qpOASES}  & $7.25/16.24$ & $0.0/0$   &   $7.71/28.26$ & $0.1/3$   &   $20.84/67.71$ & $2.8/10$ \\
\software{HPMPC}  & $15.80/17.38$ & $9.0/10$   &   $17.46/21.84$ & $9.5/11$   &   $18.95/24.43$ & $10.4/13$ \\
\software{qpDUNES}  & $3.11/10.02$ & $2.0/2$   &   $3.72/17.60$ & $2.1/4$   &   $5.24/26.90$ & $2.5/7$ \\
\software{PRESAS-SC}  & $2.30/6.39$ & $1.1/4$   &   $2.65/23.59$ & $1.2/16$   &   $4.43/43.42$ & $2.6/30$ \\
\software{PRESAS-AL}  & $2.51/7.48$ & $1.1/4$   &   $2.75/24.50$ & $1.2/16$   &   $4.87/46.28$ & $2.6/30$ \\
\software{PRESAS-PPCG}  & $1.38/2.91$ & $1.1/4$   &   $1.58/9.46$ & $1.2/16$   &   $2.50/17.73$ & $2.6/30$ \\
\software{PRESAS-PPCG-WSA}  & $1.38/2.24$ & $1.0/2$   &   $1.58/4.89$ & $1.1/6$   &   $1.81/11.41$ & $1.5/19$ \\
	\bottomrule
\end{tabular}
\vspace*{-6pt}
\end{table*}

\footnotetext{The Raspberry Pi~2 uses a BCM2836 SoC with a 900 MHz 32-bit quad-core ARM Cortex-A7 processor, with 256 KB shared L2 cache.}


\paragraph{Nonlinear MPC Simulation Results: Inverted Pendulum Swing-up}

The second numerical simulation involves a swing-up of a pendulum to its upward unstable position, using the nonlinear system dynamics and the optimal control problem formulation from~\cite{Quirynen2014a}. Both the actuated force and the cart position are constrained to remain within their respective bound values. A softened reformulation of the position constraints is used to guarantee each QP subproblem to remain feasible. The swing-up maneuver results in a relatively high number of online active-set changes with respect to the amount of state $\nx=4$ and control variables $\nU=2$. The nonlinear MPC~(NMPC) algorithm is implemented using the Gauss-Newton based RTI scheme in the \software{ACADO} code generation tool as discussed in Section~\ref{sec:NMPC_software}.

Table~\ref{tab:pend_results_v2} shows the average and worst-case computation times for each QP solution within the RTI algorithm for NMPC of the inverted pendulum, including the solvers \software{ADMM}, \software{qpOASES} and three variants of \software{PRESAS} based on single- versus double-precision arithmetics on an ARM Cortex-A7 processor. \change{The solver tolerance for \software{ADMM} was chosen to be equal to $10^{-3}$, because first-order methods are generally not efficient in computing high-accuracy solutions~\cite{Boyd2011}. Instead, a solver tolerance of $10^{-6}$ was chosen for \software{PRESAS} and the accuracy of the solutions is relatively similar to that from \software{qpOASES}.} When switching from double- to single-precision arithmetics, there can be an increased efficiency of the floating-point operations because, e.g., the compiler can more effectively use SIMD instructions. For example, this is observed for the computation time of the \software{ADMM} solver. However, the use of single-precision arithmetics may additionally lead to numerical issues that are caused by round-off errors, resulting in an overall increased number of iterations. One can observe from Table~\ref{tab:pend_results_v2} that these competing effects result in a comparable computational performance of the proposed \software{PRESAS-PPCG} solver for both single- and double-precision arithmetics, which outperforms the alternative optimal control algorithms on this particular case study. On the other hand, \software{PRESAS-SC} and \software{AL} appear to be more negatively affected by the reduced precision, while remaining competitive with the other solvers.

\begin{table}[tpb]
	\caption{NMPC for inverted pendulum swing-up~($\Tsamp = 50$~ms and $N=40$): total computation time and \# solver iterations per call, \# residual iterations per active-set change, using double- and single-precision arithmetics on an ARM Cortex-A7 processor.}
	\label{tab:pend_results_v2}
	\centering
	\setlength{\tabcolsep}{0.2em}
	\begin{tabular}{l|c|c|c}
		\toprule
		\textbf{Double precision}&{Time [ms]}&{\# sol iter}&{\# res iter} \\
		& (mean/max) & (mean/max)  	& (mean/max)\\
		\midrule
		\software{PRESAS-PPCG}&1.82/6.90&2.5/12&1.0/1\\
		\software{PRESAS-AL}&3.28/14.59&2.5/12&2.0/3\\
		\software{PRESAS-SC}&3.67/15.81&2.5/12&3.0/3\\
		\software{qpOASES}&16.96/34.86&4.2/14&-\\
		\software{ADMM}&15.87/97.57&54.2/446&-\\
		\bottomrule
	\end{tabular} \hspace{3mm}
	\begin{tabular}{l|c|c|c}
		\toprule 
		\textbf{Single precision}&{Time [ms]}&{\# sol iter}&{\# res iter} \\
		& (mean/max) & (mean/max)  	& (mean/max)\\
		\midrule
		\software{PRESAS-PPCG}&1.81/7.37&2.5/12&1.9/2\\
		\software{PRESAS-AL}&4.75/22.19&2.5/12&6.4/8\\
		\software{PRESAS-SC}&6.85/24.67&5.4/21&5.0/6\\
		\software{qpOASES}&16.37/38.86&5.2/16&-\\
		\software{ADMM}&11.68/69.61&54.5/446&-\\
		\bottomrule
	\end{tabular}
\end{table}

\subsection{Case Study~3: Hardware-in-the-Loop Simulations of NMPC for Vehicle Control} \label{sec:HIL}

The third case study illustrates the real-time computational feasibility and the memory footprint of our proposed \software{PRESAS} solver for NMPC based trajectory tracking in an autonomous driving application. The system dynamics are based on a single-track vehicle model, obtained by lumping together the left and right wheel on each axle, coupled with front and rear wheel speed dynamics and torque generation and transfer dynamics and first order dynamics for the front steering actuator. The tire friction forces are described by Pacejka's magic formula, including first order dynamics for the front and rear slip angles, and the coupling between longitudinal and lateral tire forces is based on the friction ellipse.
The trajectory tracking objective is formulated by introducing a time-dependent polynomial for each of the outputs to be tracked as an approximation of the reference motion plan. Tracking a time-dependent plan can sometimes lead to large tracking errors even when the resulting motion closely corresponds to the planned motion, e.g., caused by the vehicle slowing down or speeding up relative to the reference. Instead, the reference trajectories are parameterized with respect to a path variable, that is defined by an additional state equation and for which the change rate can be directly controlled by the NMPC controller. This results in an additional degree of freedom, such that NMPC can accelerate/decelerate along the path.
The inequality constraints in the NMPC problem formulation consist of geometric and physical limitations of the system, such as constraints on the distance of the vehicle position to the reference trajectory and on the lateral acceleration. These state-dependent constraints are formulated as soft inequality constraints, including an L1 penalty on the slack variable in the objective function. In addition, hard bounds are imposed on the steering wheel angle, steering rate, and on the wheel torques.
The resulting OCP in the NMPC controller includes $\nx=15$ differential states, $\nU=5$ control inputs and $N=20$ or $N=30$ control intervals for a $T=2$ or $T=3$~s horizon length, respectively.
The vehicle and tire friction model parameters are from~\cite{Berntorp2014}. More details on the problem formulation can be found in~\cite{Quirynen2018d,Berntorp2019a}.


We evaluate the controller in a dSPACE MicroAutoBox-II rapid prototyping unit~(RPU), which is equipped with a $900$~MHz PowerPC real-time processor~(IBM PPC 750GL) and $16$~MB of RAM. The capabilities of the RPU are slightly superior, but relatively on the same order, compared to those of current and next-generation processors for hard real-time automotive applications. We assess the real-time computational feasibility of the proposed NMPC controller based on an automatically generated C~code implementation of the RTI algorithm, using the \software{ACADO} code generation tool, in combination with different embedded QP solvers. Table~\ref{tab:timing_results} shows the average and worst-case closed-loop computation times using a sampling time of $T_{\mathrm{s}} = 100$~ms and two different values for the number of control intervals~$N_{\mathrm{MPC}}$ on the dSPACE MicroAutoBox-II. In addition, Table~\ref{tab:timing_results} presents the total memory usage for each of the software implementations of the NMPC controller. The QP solution time scales linearly with respect to the horizon length for the proposed \software{PRESAS} solvers and the total computation time remains well below the hard real-time sampling requirement of $100$~ms. Note that all three \software{PRESAS} variants in Table~\ref{tab:timing_results} are based on the \software{PPCG} solver. The results based on dense QP solvers, i.e., \software{D-PRESAS}, \software{qpOASES} and \software{ADMM}, are obtained using the dense QP interface in \software{ACADO} that includes a condensing routine to numerically eliminate the state variables from the sparse QP subproblem at each control time step. \change{As discussed in~\cite{Boyd2011}, the \software{ADMM} solver does not require the explicit use of slack variables in the OCP problem formulation to soften state-dependent inequality constraints, resulting in a reduced preparation time compared to \software{D-PRESAS} and \software{qpOASES} in Table~\ref{tab:timing_results}.} Table~\ref{tab:timing_results} shows the percentage of overruns, i.e., the time steps in which the total computation time exceeds the desired sampling time of $100$~ms. The \software{ADMM} solver occasionally does not reach the required solver tolerance within $T_{\mathrm{s}} = 100$~ms in case $N_{\mathrm{MPC}} = 20$ and the overrun percentage reaches even a value of $42$~\% in case $N_{\mathrm{MPC}} = 30$. The state of the art \software{qpOASES} solver is quite competitive with the different variants of our \software{PRESAS} solver for this NMPC simulation case study if $N_{\mathrm{MPC}} = 20$. However, for a longer control horizon, the total memory usage of the software implementation exceeds the amount of memory that is available on the dSPACE MicroAutoBox-II. On the other hand, the dense linear algebra variant of our QP solver, \software{D-PRESAS}, is more compact than \software{qpOASES}, such that this NMPC implementation fits in the dSPACE MicroAutoBox-II even for $N_{\mathrm{MPC}} = 30$. 

It can be seen from Table~\ref{tab:timing_results} that all three variants of the \software{PRESAS} solver are real-time feasible and result in a relatively low memory footprint. The sparse implementation of the \software{PRESAS-WSA} method results overall in the best computational performance, especially when looking at the worst-case computation time for a larger problem size, i.e., $N_{\mathrm{MPC}} = 30$ in Table~\ref{tab:timing_results}. \change{More specifically, one can observe that the maximum number of iterations increased from $9$ to $11$ for \software{PRESAS-WSA}, instead of the increase from $8$ to $26$ for \software{PRESAS}, when increasing the control horizon length from $N_{\mathrm{MPC}} = 20$ to $N_{\mathrm{MPC}} = 30$. The latter is a good motivation for the augmented Lagrangian type initialization in \software{PRESAS-WSA}~(see Section~\ref{sec:PRESAS_WSA}) instead of the forward simulation based initialization procedure in Section~\ref{sec:PRESAS_sim} for the \software{PRESAS} solver. Forward simulation of unstable vehicle dynamics can lead to a bad initial guess for the primal active-set strategy, which is an issue that becomes more obvious for an increasingly long prediction horizon.} For a smaller problem size, i.e., $N_{\mathrm{MPC}} = 20$, the dense implementation~(\software{D-PRESAS}) is still competitive, and its comparison with other dense solvers provides an assessment of the performance of the primal active-set strategy with the preconditioned iterative solver on its own~(without the block-structured sparsity exploitation).

\begin{table}[!ht]
	\centering
	\setlength{\tabcolsep}{0.9em}
	\begin{tabular}{l|c|c|c|c|c}
		\toprule
		$N_{\mathrm{MPC}} = 20, \; T_{\mathrm{s}} = 100$~ms & \software{PRESAS-WSA} & \software{PRESAS} & \software{D-PRESAS} & \software{qpOASES} & \software{ADMM}  	\\
		\midrule
		NMPC -- QP preparation time		&	$7.3/7.4$~ms	&	$7.3/7.4$~ms		&	$12.7/12.8$~ms			&	$12.7/12.8$~ms	&	$11.0/11.0$~ms	\\
		\midrule
		NMPC -- QP solution time		&	$12.7/22.2$~ms	&	$11.4/20.6$~ms		&	$11.3/18.0$~ms			&	$10.7/25.2$~ms	&	$43.1/137.1$~ms	\\
		NMPC -- QP solver: \# iter		&	$3/9$			&	$2/8$				&	$1/5$				&	$2/23$			&	$107/443$	\\
		\midrule
		Total time~(mean/max)		&	$20.0/29.6$~ms	&	$18.7/28.0$~ms		&	$24.0/30.8$~ms			&	$23.4/38.0$~ms	&	$54.1/148.1$~ms	\\
		Total overrun percentage		&	$0$~\%			&	$0$~\%				&	$0$~\%			&	$0$~\%			&	$7.4$~\%	\\
		Total memory usage				&	$8.56$~MB		&	$8.54$~MB			&	$8.88$~MB			&	$15.13$~MB		&	$8.46$~MB
		\vspace{5mm}\\
		\toprule
		$N_{\mathrm{MPC}} = 30, \; T_{\mathrm{s}} = 100$~ms & \software{PRESAS-WSA} & \software{PRESAS} & \software{D-PRESAS} & \software{qpOASES} & \software{ADMM}  	\\
		\midrule
		NMPC -- QP preparation time		&	$11.2/11.2$~ms	&	$11.2/11.3$~ms		&	$22.3/22.4$~ms			&	$-$			&	$18.7/18.8$~ms	\\
		\midrule
		NMPC -- QP solution time		&	$19.0/34.6$~ms	&	$18.5/71.3$~ms		&	$31.5/52.4$~ms			&	$-$			&	$75.2/197.6$~ms	\\
		NMPC -- QP solver: \# iter		&	$3/11$			&	$2/26$				&	$1/6$				&	$-$			&	$51/249$		\\
		\midrule
		Total time~(mean/max)		&	$30.2/45.8$~ms	&	$29.7/82.6$~ms		&	$53.8/74.8$~ms			&	$-$			&	$93.9/216.4$~ms	\\
		Total overrun percentage		&	$0$~\%			&	$0$~\%				&	$0$~\%			&	$-$			&	$42.4$~\%		\\
		Total memory usage				&	$9.04$~MB		&	$9.02$~MB			&	$10.1$~MB			&	$-$			&	$9.12$~MB			\\
		\bottomrule
	\end{tabular}
	\caption{Average and worst-case timing results~(ms) for a $300$~s simulation of the NMPC based vehicle control algorithm using a sampling time $T_{\mathrm{s}} = 100$~ms and $N_{\mathrm{MPC}} = 20$ or $N_{\mathrm{MPC}} = 30$ on dSPACE MicroAutoBox-II. All three \software{PRESAS} variants are based on the \software{PPCG} solver and the results for the dense QP solvers, \software{D-PRESAS}, \software{qpOASES} and \software{ADMM}, are obtained using condensing.}
	\label{tab:timing_results}
\end{table}


\section{Experimental Results for NMPC of Autonomous Scaled Vehicles} \label{sec:hamsters}

In this section, we present experimental results for an NMPC implementation, using the proposed \software{PRESAS} solver, of a real-time trajectory tracking controller for an autonomous driving system, implemented on a testbench of small-scale vehicles~\cite{Berntorp2018,Berntorp2019}. Such a platform, based on the Hamster robots as illustrated in Fig.~\ref{fig:hamster}, can be used for rapid prototyping and evaluation of a complete control and estimation software stack \change{before deploying to full-scale vehicle experiments, e.g., see~\cite{Berntorp2019b,Uno2019}.} Next, we briefly describe the test setup and the NMPC problem formulation, before presenting the experimental results in Section~\ref{sec:exp_results}.

\subsection{Scaled Vehicle Experimental Platform}

\begin{figure*}
	\centering
	\begin{minipage}{.3\linewidth}%
				\vspace{6mm}
		\begin{subfigure}[b]{\textwidth}
			\centering
			\includegraphics[width=0.8\textwidth,trim={0 15cm 0 15cm},clip]{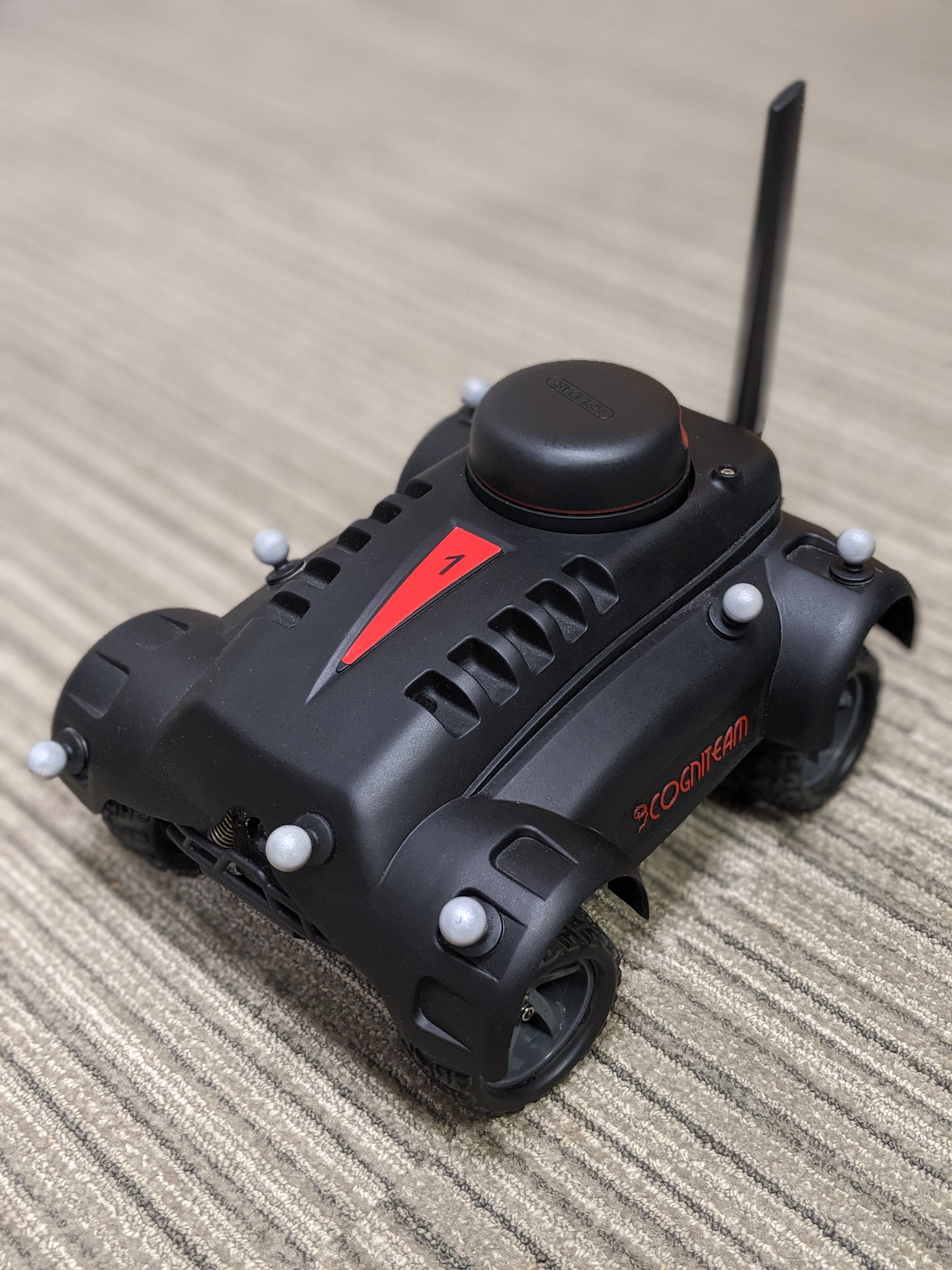}
			\captionsetup{font=normalsize,labelfont={sf}}
			\caption[]%
			{\normalsize The Ackermann-steered Hamster mobile robot used in the experiments.}    
			\label{fig:hamster}
		\end{subfigure}
		\vskip\baselineskip\vspace{2mm}
		\begin{subfigure}[b]{\textwidth}   
			\centering 
			\includegraphics[width=0.8\textwidth,trim={0 15cm 0 15cm},clip]{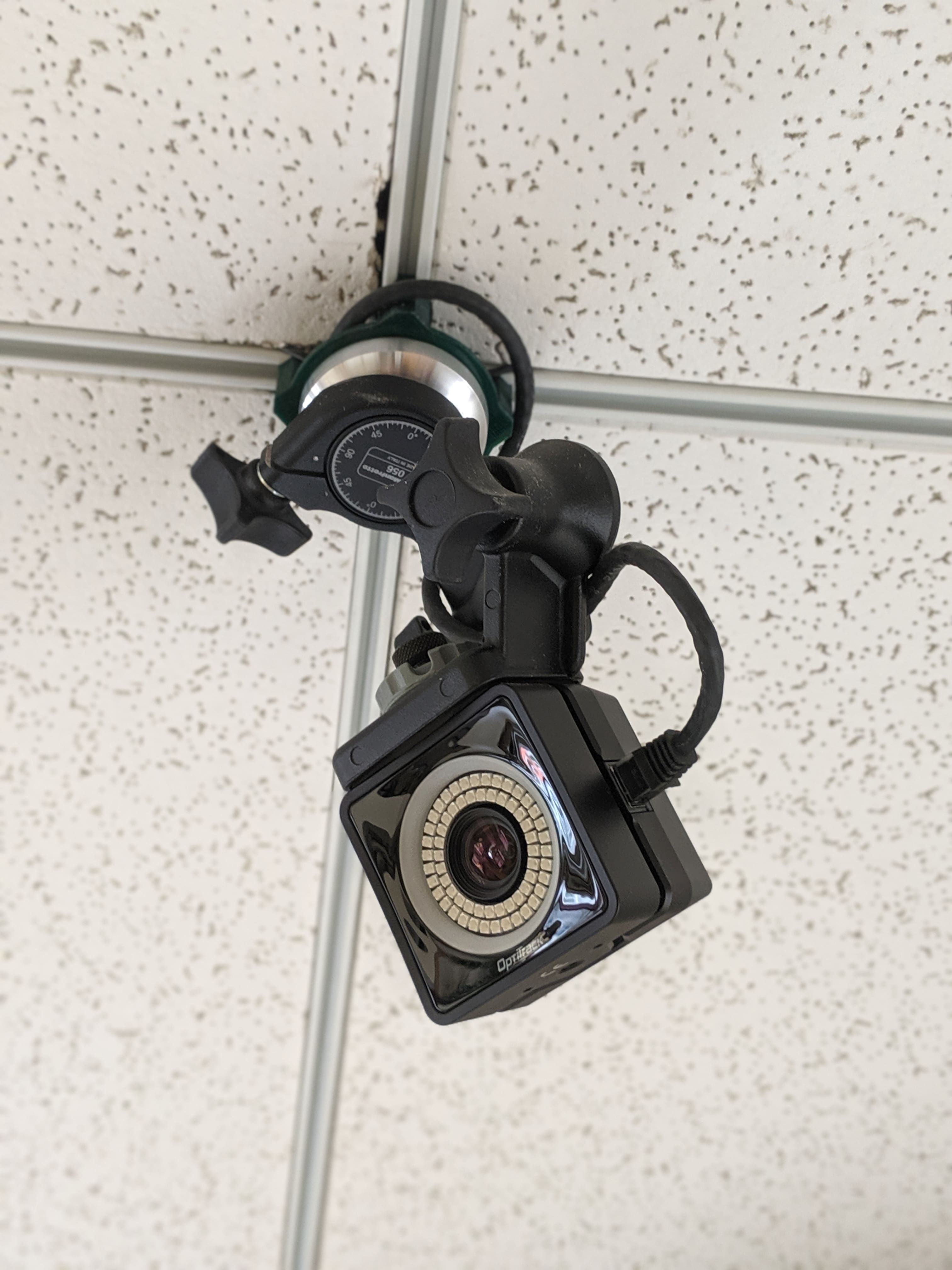}
			\captionsetup{font=normalsize,labelfont={sf}}
			\caption[]%
			{\normalsize A camera from the OptiTrack motion capture system.}    
			\label{fig:optitrack}
		\end{subfigure}
	\end{minipage}
	\begin{minipage}{.6\linewidth}%
		\begin{subfigure}[b]{\textwidth}  
			\centering 
			\includegraphics[trim={0 3.5cm 0 0},clip,width=1.1\textwidth]{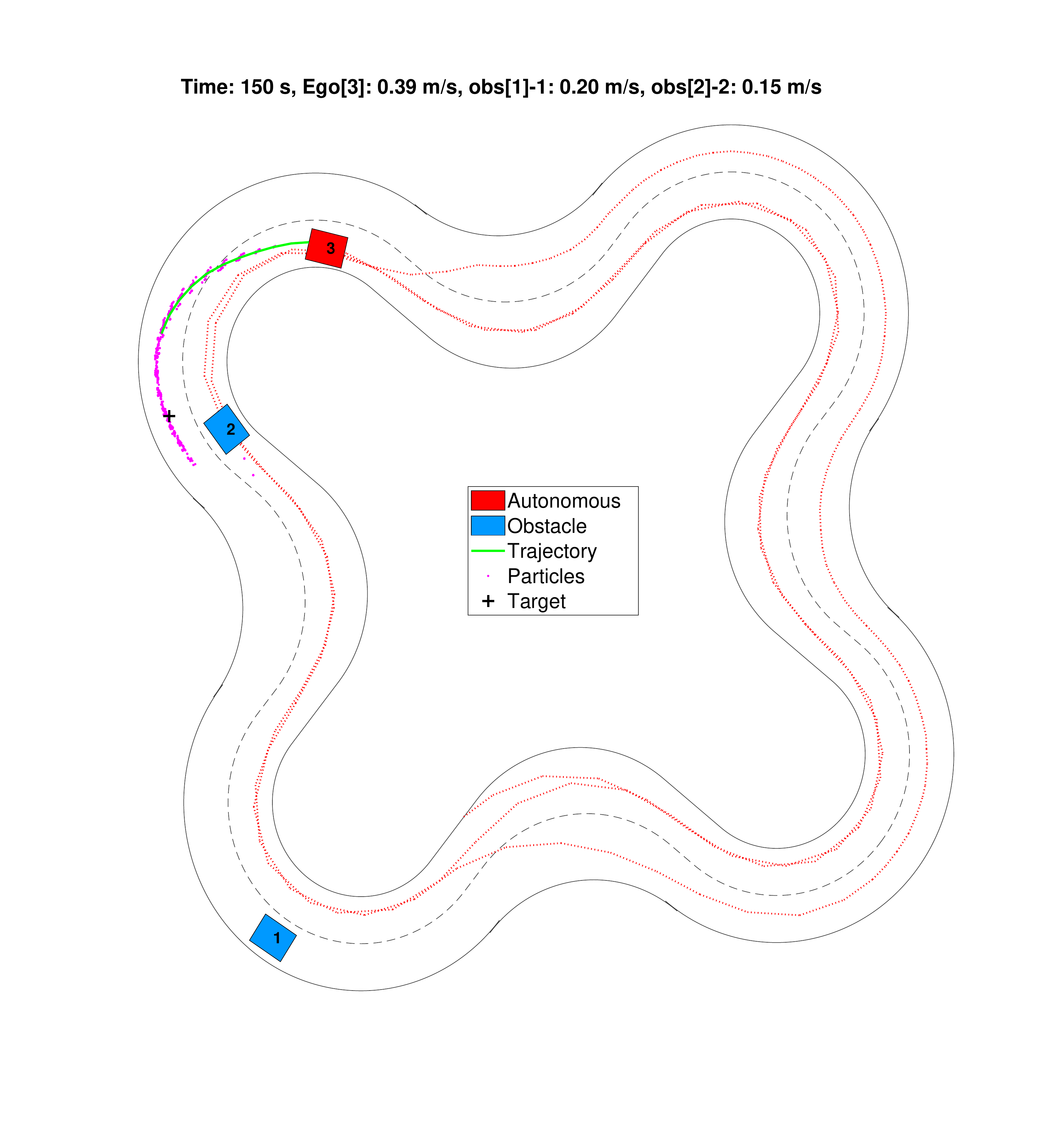}
			\captionsetup{font=normalsize,labelfont={sf}}
			\caption[]%
			{\normalsize Current pose of moving obstacles~(in blue) and autonomous vehicle~(in red), closed-loop trajectory~(dotted red) and predicted motion plan~(in green).}
			\label{fig:track}
		\end{subfigure}
	\end{minipage}
	\caption[ ]{Illustration of our experimental testbench that consists of small-scale test vehicles~(a), with on-board sensors, and an OptiTrack motion capture system~(b). Figure~(c) shows the track and a snapshot of the position of the autonomous vehicle and of the two obstacle vehicles, as well as an illustration of the closed-loop trajectory at the end of the experiment.}
	\label{fig:hamster_setup}
\end{figure*}

The Hamster is a $25\times20$~cm mobile robot for research and prototype development, see Fig.~\ref{fig:hamster}. It is equipped with scaled versions of sensors commonly available on full-scale research vehicles, such as a $6$~m range mechanically rotating $360$~deg Lidar, an inertial measurement unit~(IMU), GPS receiver, HD camera, and motor encoders. It uses two Raspberry Pi~3 computing platforms, each with an ARM Cortex-A53 processor. The robot has Ackermann steering and therefore it is kinematically equivalent to a full-scale vehicle, and its dynamic behavior resembles that of a regular vehicle. The Hamster has low-level dedicated hardware for power distribution and monitoring. At the lowest layer, the robot is controlled by setting the desired front wheel steering angle and longitudinal velocity.  In addition, the Hamster has built-in mapping and localization capabilities, and object detection and tracking can be done with the onboard Lidar and/or camera. Hence, the platform is a good test setup for verifying the dynamic feasibility and performance of vehicle control and estimation software in a realistic setting, with a sensor setup similar to the one expected in full-scale autonomous vehicles. The Hamster communicates and connects to external algorithms using the robot operating system (ROS). 
In order to be able to accurately evaluate the control and estimation algorithms in terms of performance in a controlled environment, we use an OptiTrack motion capture system. The OptiTrack system is a flexible camera-based (see Fig.~\ref{fig:optitrack}) six degrees-of-freedom tracking system that can be used for tracking drones, ground vehicles, and industrial robots. Depending on the environment and quality of the calibration, the system can track the position of the Hamster robot within $0.9$~mm accuracy and with a rotational error of less than $3$~deg. The OptiTrack is connected to the same ROS network using the common VRPN protocol.


\subsection{Vehicle Dynamics and NMPC Formulation}

Due to the testing environment, we consider a kinematic single-track model in which the two wheels on each axle are lumped together. Although a dynamic vehicle model based on force balances is generally more accurate than a kinematic model, the latter may be still sufficient for several driving conditions~\cite{Carvalho2015}, and model errors are corrected for by the feedback control. On the other hand, a dynamic model depends on more parameters, such as the wheel radii, tire stiffness, and vehicle mass and inertia, which typically are unknown/uncertain and may be difficult, or at least tedious, to estimate. Here, for simplicity, we use the kinematic single-track model
\begin{equation}\label{eq:vehmod}
\dot x = \begin{bmatrix}
\dot{p}_{\mathrm{X}} \\
\dot{p}_{\mathrm{Y}} \\[2pt]
\dot{\psi} \\ 
\dot{\delta}_{\mathrm{f}} \\
\dot{v}_{\mathrm{x}} \\
\dot{\delta}
\end{bmatrix} = \begin{bmatrix}
v_{\mathrm{x}}\cos(\psi + \beta)\\[2pt]
v_{\mathrm{x}}\sin(\psi + \beta) \\[2pt] 
\frac{v_{\mathrm{x}}}{L} \tan(\delta_{\mathrm{f}}) \cos(\beta) \\
1/t_{\mathrm{d}}(\delta + \delta_0 - \delta_{\mathrm{f}}) \\
u_1 \\
u_2
\end{bmatrix},
\end{equation}
where $p_{\mathrm{X}},p_{\mathrm{Y}}$ is the longitudinal and lateral position in the world frame, $\psi$ is the heading angle and $\dot \psi$ the heading rate of the vehicle, $v_{\mathrm{x}}$ is the longitudinal velocity of the vehicle, $\delta$ and $\delta_{\mathrm{f}}$ are, respectively, the desired and actual front wheel steering angle, $L:=l_f+l_r$ is the wheel base, and $\beta := \arctan{ ({l_r \tan(\delta_{\mathrm{f}})}/{L})}$ is the body-slip angle. A first order front steering equation is included in Eq.~\eqref{eq:vehmod} which models that, due to the steering mechanism and other platform effects, the wheel angle response is not immediate. In addition, an input additive offset value $\delta_0$ can be estimated online. The inputs $u_1,u_2$ are the acceleration and the steering rate, respectively. This choice allows to provide smooth velocity and steering profiles and to constrain the allowed change rate of the vehicle velocity and front steering wheel angle.

Similar to the vehicle control hardware-in-the-loop simulation results in Section~\ref{sec:HIL}, the NMPC trajectory tracking objective is formulated based on a smooth approximation of the reference motion plan at each control time step, parameterized with respect to a time-dependent path variable $\tau$ for which the change rate $\frac{\dd \tau}{\dd t} = \dot{\tau}$ is an additional control input. In continuous time, the NMPC cost function therefore consists of the following terms
\begin{equation}
\int_{0}^{T} \! \left( \Vert F^{\mathrm{ref}}(x(t)) - y^\mathrm{ref}(\tau,\polpar) \Vert_{W}^2 \;+\; 
\Vert x(t) \Vert_{Q}^2 \;+\; \Vert u(t) \Vert_{R}^2 \;+\; r_{\mathrm{s}}\, s(t) \right) \, \dd t,
\end{equation}
including a term for tracking the reference motion plan and two regularization terms for penalizing state and control variables, where $Q \in \R^{\nx \times \nx}$ and $R \in \R^{\nU \times \nU}$ are the corresponding weighting matrices. The definition of a positive slack variable $s(t) \ge 0$ allows for an easy implementation of an exact L1 soft constraint penalty in the objective, i.e., $| s(t) | = s(t)$. The reference motion plan is approximated by a smooth function $y^\mathrm{ref}(\tau,\polpar)$ that depends on the path variable $\tau$ as well as on additional parameters $\polpar$. We define the tracking function in the objective based on a polynomial approximation of the reference trajectories as
\begin{equation}
F^{\mathrm{ref}}(x(t)) - y^\mathrm{ref}(\tau,\polpar) = \begin{bmatrix} e_{Y}(t,\tau,\polpar) 
\\ p_{X}(t)-p_{X}^\mathrm{ref}(\tau,\polpar_X)  \\ p_{Y}(t)-p_{Y}^\mathrm{ref}(\tau,\polpar_Y) \\ \psi(t)-\psi^\mathrm{ref}(\tau,\polpar_{\psi}) \\ 
v_{\mathrm{x}}(t)-v_{\mathrm{x}}^\mathrm{ref}(\tau,\polpar_v) \end{bmatrix} = \begin{bmatrix} e_{Y}(t,\tau,\polpar_X,\polpar_Y,\polpar_\psi) \\ p_{X}(t) - \sum_{i=0}^{n_X}\polpar_X^i \tau^i  \\ p_{Y}(t) - \sum_{i=0}^{n_Y}\polpar_Y^i \tau^i \\ \psi(t) - \sum_{i=0}^{n_\psi}\polpar_\psi^i \tau^i \\ v_{\mathrm{x}}(t) - 
\sum_{i=0}^{n_v}\polpar_v^i \tau^i \end{bmatrix},
\end{equation}
where the path error $e_{Y}(\cdot) = \text{cos}(\psi^\mathrm{ref}(\tau,\polpar_\psi))\left(p_{Y}-p_{Y}^\mathrm{ref}(\tau,\polpar_Y)\right) - \text{sin}(\psi^\mathrm{ref}(\tau,\polpar_\psi))\left(p_{X}-p_{X}^\mathrm{ref}(\tau,\polpar_X)\right)$ is defined as the orthogonal distance to the parameterized reference trajectory. 
The inequality constraints in the NMPC problem formulation include hard bounds on the control inputs and soft constraints for limiting the distance to the parameterized reference trajectory, the vehicle velocity and the front steering wheel angle:
\begin{subequations}
	\begin{alignat}{5}
	1 - \Delta\bar{\tau} &\le  \dot{\tau} \le 1 + \Delta\bar{\tau}, \qquad
	&-\bar{\dot{\delta}} &\le \; \dot{\delta} \le 
	\bar{\dot{\delta}}, \qquad
	&-\bar{\dot{v}}_{\mathrm{x}} &\le \; \dot{v}_{\mathrm{x}} \le 
	\bar{\dot{v}}_{\mathrm{x}}, \qquad &0 &\,\le \, s, \\
	-\bar{e}_Y &\le \; e_{Y} + s, \qquad
	& -\bar{\delta}_f &\le \; \delta_f + s, \qquad
	& -\bar{v}_{\mathrm{x}} &\le \; v_{\mathrm{x}} + s, \\
	e_{Y} - s &\le \bar{e}_Y, \qquad
	& \delta_f - s &\le \bar{\delta}_f,\qquad
	& v_{\mathrm{x}} - s &\le \bar{v}_{\mathrm{x}}.
	\end{alignat}\label{eq:softConstr}
\end{subequations}
The resulting OCP includes $\nx=7$ differential states, $\nU=4$ control inputs and $N=60$ control intervals with a sampling time of $T_{\mathrm{s}} = 25$~ms over a $T=1.5$~s horizon length. The NMPC controller is implemented with a sampling frequency of $40$~Hz, using the direct multiple shooting discretization method in combination with the RTI scheme in the \software{ACADO} code generation tool and our proposed \software{PRESAS} optimization algorithm.

\subsection{Experimental Results and Solver Performance}
\label{sec:exp_results}

We use three Hamsters in the experimental validation, one acts as the controlled ego vehicle~(EV) and two Hamsters act as dynamic obstacles. The overall objective is to avoid the obstacles while circulating a two-lane closed circuit in the counter-clockwise direction, see Fig.~\ref{fig:track}, with the inner lane as preferred lane. The obstacle vehicles~(OVs) are commanded to track the middle of either of the lanes, where the preferred lane and vehicle velocities can be changed throughout the experiments. The OVs use PID controllers for tracking the desired lane and velocity. The current position of each OV is obtained from the OptiTrack while its velocity is estimated. 
For the EV, the desired path trajectory, front steering angle, and vehicle velocity are computed by the motion planner and sent to the NMPC controller that tracks the trajectory in real time at a sampling frequency of $40$~Hz. For this purpose, a particle filtering based motion planner is used, as presented recently in~\cite{Berntorp2019}. The planning algorithm typically runs at a relatively low frequency of $1$-$2$~Hz and therefore relies on the NMPC controller to make the vehicle accurately execute the desired motion. The current state of the vehicle, as well as disturbances such as the steering offset value, are estimated online using an extended Kalman filter based on IMU and OptiTrack measurements.

\begin{figure}[tpb]
	\centerline{\hbox{
			\includegraphics[width=0.95\textwidth]{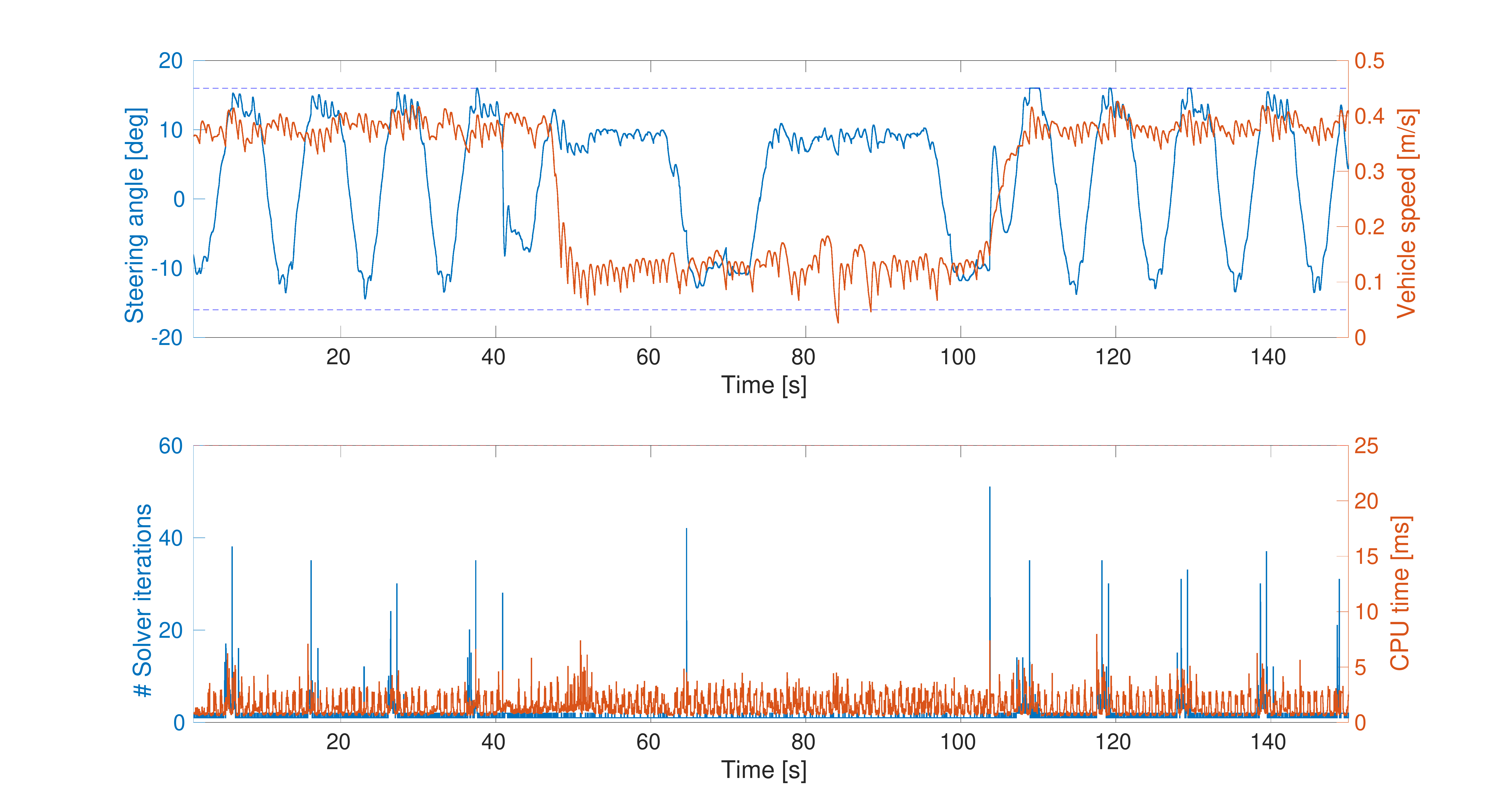}}}
	\vspace*{-8pt}
	\caption{Closed-loop results from Hamster experiments~(see Figure~\ref{fig:hamster_setup}): commanded front steering wheel angle $\delta$, vehicle speed $v_{\mathrm{x}}$, number of \software{PRESAS} solver iterations and computation time at each control time step.}
	\vspace*{-8pt}
	\label{fig:hamster_results}
\end{figure}

Figure~\ref{fig:hamster_results} shows closed-loop results from $150$~s of Hamster experiments on the test track in Figure~\ref{fig:track}. The upper plot shows time trajectories of the control inputs that were sent by the NMPC controller to the actuators of the EV. The target velocity for the EV was chosen to be $0.4$~m/s. It can be observed that the EV needed to slow down to a velocity of about $0.15$~m/s between $50$ and $100$~s, since both lanes were blocked by slow moving OVs during this time period. The lower plot of Figure~\ref{fig:hamster_results} illustrates the computational performance of our proposed \software{PRESAS-PPCG-WSA} implementation in Algorithm~\ref{alg:aug_lag}. Due to the warm starting of the solver, the number of active-set changes remains typically low at each control time step. Occasionally, due to a bad initial guess, the number of solver iterations can jump up to a relatively high number. But even the worst-case performance results in a total computation time that is well below the desired sampling time of $T_{\mathrm{s}}=25$~ms in these experiments.


\section{Conclusions} \label{sec:concl}

This paper presented different variants of a tailored implementation of a primal active-set strategy for solving optimal control structured quadratic programming problems as they typically arise in both linear and nonlinear real-time MPC applications. The proposed \software{PRESAS} solver is based on either a sparse block-diagonal preconditioner for the preconditioned minimal residual~(PMINRES) method or a block-sparse constraint preconditioner for the projected preconditioned conjugate gradient~(PPCG) method within a primal active-set optimization algorithm. In addition, we described different warm-started initialization procedures to compute a primal feasible starting point in a real-time predictive controller. Based on an efficient and self-contained C~code implementation, the computational performance and numerical reliability of \software{PRESAS} is illustrated against other state of the art QP solvers. We also showed that the \software{PRESAS} QP~solver is real-time feasible on a dSPACE MicroAutoBox-II rapid prototyping unit for vehicle control applications. Finally, we presented experimental results for a real-time NMPC implementation based on \software{PRESAS} on a testbench of small-scale autonomous vehicles.

\bibliography{syscop}

\end{document}